\documentclass{article}
\usepackage[utf8]{inputenc}
\usepackage{lipsum}
\usepackage{amsfonts}
\usepackage{graphicx}
\usepackage{epstopdf}
\usepackage{algorithmic}
\usepackage{amsmath}
\usepackage{amssymb}
\title{LQG control over unreliable channels-Full Proof}
\author{Fredrik Bengtsson and Torsten Wik}

\begin{document}

\maketitle

\section{Introduction}
The use of wireless communication between sensors, actuators and controllers
can create large savings by avoiding costly wiring and hardware. Depending
on locations and external conditions, however, it can be difficult
to ensure reliable communication. This raises numerous issues for control and estimation, as discussed in \cite{hespanha2007survey}. 

LQG control (optimized state feedback control for linear systems with
Gaussian disturbance and quadratic cost) is a well established method
applicable to MIMO systems. It was introduced already in the sixties
and remains one of the most implemented type of controllers. Consequently,
it is important to examine how to determine an LQG controller when
communication channels are unreliable.

When dealing with LQ control over unreliable links there are several different aspects to consider (see Table \ref{aspects}), each implying different solutions. 
\begin{table}
\centering
\tiny
\caption{LQG control over unreliable links - problem formulations. Here we will examine finite-horizon control in the TCP like case using hold-input with both packet losses and delays.}

\label{aspects}
\begin{tabular}{|l|c|c|}
\hline 
\textbf{Type} & Packet loss & Packet delays\tabularnewline
\hline 
\textbf{Control} & Finite horizon & Infinite horizon\tabularnewline
\hline 
\textbf{Strategy} & Zero-input & Hold input\tabularnewline
\hline 
\textbf{Signal} & No acknowledgements(UDP) & Acknowledgements(TCP)\tabularnewline
\hline 
\end{tabular}

\end{table}

\normalsize
First of all there are two basic ways a channel can be unreliable, one being \textit{packet
drops} and the other \textit{packet delays}. LQG control in the case of packet
losses has been examined in \cite{Imer20061429,1184678,4118476,5717722}. Then the cost can be over a finite or infinite time horizon.
Infinite horizon LQG control for the delayed case has been examined
in \cite{lincoln2000optimal} and \cite{shousong2003stochastic}, where it is shown that the optimal
solution will not only depend on the states but also on the previous
control signals. In \cite{wang2018optimal} and \cite{lincoln2000optimal} infinite horizon control for the case of a system with both random time delays and packet dropouts is investigated. However, no explicit way to derive the solution is
presented. In \cite{liang2016optimal} and \cite{ma2017optimal}  finite horizon LQ control is examined for the case of packet losses and constant delays. 

When there is an unreliable channel between the controller and the
actuator this means that the latest signal sent might not yet have
arrived when the actuator needs to execute a new control action. When this occurs there are two basic strategies the actuator
can adopt \cite{6760935,4908929}. One is to apply the last input
received until a more recent one arrives, which is known as \textit{hold-input}.
In the other one, the input is set to zero if the latest control signal
sent is delayed or lost, which is known as \textit{zero-input}. Zero-input
treats delayed and lost signals the same way and zero-input optimizing
control for a system with a lossy channel has been derived
in \cite{Imer20061429}. The focus here though is on designing a hold-input
control strategy for a system subject to a random unbounded delay
and to packet losses.

Finally, there are two basic types of unreliable communication links. In one
the sender does not know if the sent packet has arrived, which is
known as the \textit{UDP-like case}. In the other one there is a system of
acknowledgment that ensures that the sender knows if the sent packet
has arrived. This is referred to as the \textit{TCP-like case}. In this article we will focus on the TCP-like case. Furthermore we will assume the acknowledgments arrive without losses or delays.  In previous work LQG control
for this case was examined  \cite{bengtsson2016lqg}, but
only for a specific probability function for the delay. Moreover,
 the knowledge whether the sent packet
has arrived or not was then only used to facilitate the estimation of the
states. Here, this knowledge is also used to optimize the control
scheme, yielding a more optimal solution. This optimality, though,
comes at the price of increased complexity of the solution and increased
computational cost.

As mentioned, we will examine a system that can be affected both
by random unbounded delays and packet losses, and derive an explicit solution for
finite horizon optimal LQG control. The assumption will be that the probability functions of the delays and packet losses are known. If the probability functions are unknown Q-learning can be used to derive solutions as discussed in \cite{xu2012stochastic}. The focus  is to present a detailed and complete derivation of a implementable solution to the discussed optimal control problem.

 \section{Problem formulation}
 \label{sec:problem}

The plant considered is assumed to be an LTI system on the form

\begin{equation}
x_{k+1}=Ax_{k}+Bu_{k}+w_{k},\label{eq:3}
\end{equation}
where $x\in R^{n}$ is the state vector, $w\in R^{n}$ is white Gaussian
noise, and $u_{k}\in R^{m}$ is the control signal applied by the
actuator at time $k$ (for hold-input control this will be the latest
control signal that has arrived). As the noise is white, and the LQG cost criterion is with respect to expectation, the noise will not
impact the optimal control scheme and therefore we will disregard
it from now on. 

The communication between the controller and the actuator is subject
to a random delay. This delay is assumed to be described by a known probability
function $p(d)$, where $d$ is the number of samples the
packet is delayed. Note that there are no requirements on $p(d)$;
it can be of both finite or infinite length. Furthermore, $p(d)$
can be used to describe systems with both delays and packet losses.
In the latter case $\left\{ p(d)\right\} _{d=0}^{\infty}$ will simply sum up to less than one.

The delays between consecutive steps are assumed to be independent.
From this assumption the probability that the latest control signal
that has arrived is the signal sent $i$ time units before can be
derived as

\begin{equation}
p_{d}(i)=P(i)\prod_{j=0}^{i-1}\bar{P}(j), \label{pddef2}
\end{equation}
where $P(i)=\sum_{k=0}^{i}p(k)$ denotes the cumulative probability and $\bar{P}(i)$ denotes
the complementary probability to $P(i)$, i.e. $\bar{P}(i)$=1-$P(i)$, which is the probability that none of the latest $i$ signals have arrived.



Another way to express  $p_{d}(i)$ is

\begin{equation}
p_{d}(i)=P(i)\bar{P}_{d}(i-1),
\label{pddef3}
\end{equation}
where $\bar{P}_{d}(i-1)$ is the complementary probability of the cumulative probability $P_{d}(i-1)=\sum_{k=0}^{i-1}p_{d}(k)$.

Now, the goal is to design a controller that determines the control
signals $v$, sent from the controller, which minimizes the quadratic
criterion

\begin{equation}
J_{N}=E\left[\sum_{i=0}^{N}u_{i}^{T}Ru_{i}+\sum_{i=0}^{N}x_{i}^{T}Qx_{i}+x_{N+1}^{T}S_{N+1}x_{N+1}\right],\label{eq:6}
\end{equation}
where $R$ is a positive definite symmetric matrix and $Q$ and $S_{N+1}$
are positive semi-definite symmetric matrices. 

As previously mentioned, we will examine the TCP-like case, where the
controller knows if a signal has reached the actuator or not. To handle this we use a variable size controller state $\zeta_{k}$  holding all issued control signals since the last one acknowledged, defined by
\[
\zeta_{k}=\left[\begin{array}{c}
v_{k}\\
\vdots\\
v_{\tau_{k}}\\
x_{k+1}
\end{array}\right],
\]
where $u_{k}=v_{\tau_{k}}$ is the control signal sent at time $\tau_k$. The update $\zeta_{k+1}$ of the state $\zeta_k$ follows from the update of $\tau_k$, i.e.
\begin{equation}
    \tau_{k+1} = \begin{cases}
                \tau, & \text{if $v_{\tau}, \,\tau>\tau_k,$ is the most recently acknowledged signal} \\
                \tau_k, & \text{if no, more recent, $v_{}$ is acknowledged}
            \end{cases}
\end{equation}
Note that in the update of $\tau_k$, obsolete acknowledgements are always discarded. 
\section{Full proof}\label{sec:fullproof}

In this section we will present the solution and describe how it is derived. As for
the derivation of the standard LQG solution we will use dynamic programming.
This means that we will start by finding the last optimal control
signal $v_{N}$ that minimizes the cost function (\ref{eq:6}), expressed
in terms of states and previous control signals available at that
time i.e. $\zeta_{N-1}$. For the remaining cost we will then find the control signal
$v_{N-1}$ that minimizes this cost expressed in signals available
at that time, $\zeta_{N-2}$. After this, $v_{N-2}$ is found to minimize the now
remaining cost. Repeating this once more reveals a pattern of induction
such that all the remaining $v_{k}$ can be calculated. 

We start by noting that 

\[
E\left[u_{i}^{T}Ru_{i}\right]=\sum_{j=0}^{i}p_{d}(i-j)v_{j}Rv_{j}.
\]

From this, the first term in (\ref{eq:6}) becomes

\begin{align*}
E\left[\sum_{i=0}^{N}u_{i}^{T}Ru_{i}\right]= & E\left[\sum_{i=0}^{N}\sum_{j=0}^{i}p_{d}(i-j)v_{j}^{T}Rv_{j}\right]\\
= & E\left[\sum_{i=0}^{N}\sum_{j=0}^{N-i}p_{d}(j)v_{i}^{T}Rv_{i}\right]\\
= & E\left[\sum_{i=0}^{N}P_{d}(N-i)v_{i}^{T}Rv_{i}\right]\\
= & E\left[\sum_{i=0}^{N}v_{i}^{T}R_{i}v_{i}\right],
\end{align*}
where

\[
R_{i}\triangleq P_{d}(N-i)R.
\]

The term $x_{N+1}^{T}S_{N+1}x_{N+1}$ we wish to express in terms of signals available
at time $N$, so using the state equation (\ref{eq:3}),

\begin{align*}
E\left[x_{N+1}^{T}S_{N+1}x_{N+1}\right]=&  E\left[u_{N}^{T}B^{T}S_{N+1}Bu_{N}\right.\\
&+2u_{N}^{T}B^{T}S_{N+1}Ax_{N}\\
 & \left.+x_{N}^{T}A^{T}S_{N+1}Ax_{N}\right]
\end{align*}
where
\begin{align}
E\left[u_{N}^{T}B^{T}S_{N+1}Bu_{N}\right]=&\sum_{i=0}^{N}p_{d}(N-i)v_{i}^{T}B^{T}S_{N+1}Bv_{i}\nonumber \\
E\left[u_{N}^{T}B^{T}S_{N+1}Ax_{N}\right]= & E\left[p_{d}(0)v_{N}^{T}B^{T}S_{N+1}Ax_{N}\right.\nonumber \\
 & \left.+\bar{P}_{d}(0)u_{N|u_{N}\neq v_{N}}^{T}B^{T}S_{N+1}Ax_{N}\right]\label{eq:newsak-1}
\end{align}
and  $E\left[u_{N|u_{N}\neq v_{N}}\right]$ is the expected value
of the actuated control signal at time $N$ given that the control
signal that was sent at time $N$ has not yet arrived. Thus

\begin{align}
E\left[x_{N+1}^{T}S_{N+1}x_{N+1}\right]= & E\left[\sum_{i=0}^{N}p_{d}(N-i)v_{i}^{T}B^{T}S_{N+1}Bv_{i}\right.\nonumber \\
 & +2p_{d}(0)v_{N}^{T}B^{T}S_{N+1}Ax_{N}\nonumber \\
 & +2\bar{P}_{d}(0)u_{N|u_{N}\neq v_{N}}^{T}B^{T}S_{N+1}Ax_{N}\nonumber \\
 & +x_{N}^{T}A^{T}S_{N+1}Ax_{N}\bigg].\label{eq:xSx}
\end{align}

Now, define
\begin{align}
T(0,b)\triangleq & \sum_{i=0}^{N+1-b}v_{i}^{T}R_{i}v_{i}+\sum_{i=0}^{N+1-b}p_{d}(N-i)v_{i}^{T}B^{T}S_{N+1}Bv_{i}\nonumber \\
 & +x_{N}^{T}A^{T}S_{N+1}Ax_{N},\label{eq:TN}
\end{align}
where $b$ is a counter. Increasing it by one removes the latest control signal from the expression. The first argument (0), is used in conjunction with $b$ to specify which is the latest control signal contained in the expression. For example, in this case the latest control signal $T(0,b)$ contains is at $N+1-b-0$. The criterion (\ref{eq:6}) can then be written as
\begin{align}
J_{N}= & E\bigg[T(0,1)+2p_{d}(0)v_{N}^{T}B^{T}S_{N+1}Ax_{N}\nonumber \\
 & \left.+2\bar{P}_{d}(0)u_{N|u_{N}\neq v_{N}}^{T}B^{T}S_{N+1}Ax_{N}+\sum_{i=0}^{N}x_{i}^{T}Qx_{i}\right].\label{eq:7-1-1}
\end{align}

To extract the parts of the cost that depend on the latest control
signal we first note that

\[
T(0,b)=T(0,b+1)+v_{N+1-b}^{T}T_{c}(0,b)v_{N+1-b},
\]
where

\begin{align}
T_{c}(0,b)\triangleq R_{N+1-b}+p_{d}(b-1)B^{T}S_{N+1}B. \label{TC0}
\end{align}

The cost (\ref{eq:7-1-1}) can then be expressed as

\begin{align}
J_{N}= & E\bigg[v_{N}^{T}A_{11}(N)v_{N}+2v_{N}^{T}A_{12}(N)x_{N}\nonumber \\
 & +2\bar{P}_{d}(0)u_{N|u_{N}\neq v_{N}}^{T}B^{T}S_{N+1}Ax_{N}\nonumber \\
 & \left.+T(0,2)+\sum_{i=0}^{N}x_{i}^{T}Qx_{i}\right],\label{eq:7-1}
\end{align}
where

\begin{eqnarray*}
A_{11}(N) & = & T_{c}(0,1)\\
A_{12}(N) & = & p_{d}(0)B^{T}S_{N+1}A
\end{eqnarray*}

Only the first two terms in (\ref{eq:7-1}) depend on $v_{N}$ and as $A_{11}$ is positive definite
the optimal $v_{N}$ that minimizes (\ref{eq:7-1}) is

\begin{equation}
v_{N}=-A_{11}^{-1}(N)A_{12}(N)x_{N},\label{eq:optimal}
\end{equation}
which results in a minimum cost 

\begin{align}
J_{N}^{*}= & E\bigg[-x_{N}^{T}(A_{12}^{T}(N)A_{11}^{-1}(N)A_{12}(N))x_{N}+T(0,2)\nonumber \\
 & \left.+2\bar{P}_{d}(0)u_{N|u_{N}\neq v_{N}}^{T}B^{T}S_{N+1}Ax_{N}+\sum_{i=0}^{N}x_{i}^{T}Qx_{i}\right].\label{eq:costleft2}
\end{align}

The next step is to find the control signal $v_{N-1}$ that minimizes (\ref{eq:costleft2}). Let us start by examining

\begin{equation}
E\left[T(0,2)-x_{N}^{T}(A_{12}^{T}(N)A_{11}^{-1}(N)A_{12}(N))x_{N}+\sum_{i=0}^{N}x_{i}^{T}Qx_{i}\right],\label{eq:TNsak}
\end{equation}
from which we want to extract the part that depends on the latest
state $x_{N}$, such that we can use the state equation to eliminate
$x_{N}$. Therefore, we split the function $T$, defined by (\ref{eq:TN}),
as

\begin{equation}
T(0,2)=T_{noX}(0,2)+x_{N}^{T}T_{X}(0)x_{N},\label{TNxremove}
\end{equation}
where

\[
T_{X}(0)=A^{T}S_{N+1}A
\]
and

\[
T_{noX}(0,b)=\sum_{i=0}^{N+1-b}v_{i}^{T}R_{i}v_{i}+\sum_{i=0}^{N+1-b}p_{d}(N-i)v_{i}^{T}B^{T}S_{N+1}Bv_{i}.
\]
The expression (\ref{eq:TNsak}) then becomes

\[
E\left[T_{noX}(0,2)+x_{N}S_{N}x_{N}+\sum_{i=0}^{N-1}x_{i}^{T}Qx_{i}\right],
\]
where

\[
S_{N}=T_{X}(0)-A^{T}_{12}(N)A_{11}^{-1}(N)A_{12}(N)+Q.
\]

By expanding $x_{N}S_{N}x_{N}$ using the state equation, in the same way as was done when deriving (\ref{eq:xSx}) we get

\begin{eqnarray*}
E\left[x_{N}S_{N}x_{N}\right] & = & E\left[2p_{d}(0)v_{N-1}^{T}B^{T}S_{N}Ax_{N-1}+2\bar{P}_{d}(0)u_{N-1|u_{N-1}\neq v_{N-1}}^{T}B^{T}S_{N}Ax_{N-1}\right.\\
 &  & \left.+\sum_{i=0}^{N-1}p_{d}(N-1-i)v_{i}^{T}B^{T}S_{N}Bv_{i}+x_{N-1}^{T}A^{T}S_{N}Ax_{N-1}\right]
\end{eqnarray*}
for $k>0$, we define


\begin{eqnarray}
    T(k,b) &= &T_{noX}(k-1,b+1)+\sum_{i=0}^{N-k+1-b}p_{d}(N-k-i)v_{i}^{T}B^{T}S_{N-k+1}Bv_{i}\nonumber\\
    & &+x_{N-k}^{T}A^{T}S_{N-k+1}Ax_{N-k}\label{T_kb}
\end{eqnarray}
 and thus (\ref{eq:TNsak}) becomes

\begin{align*}
 & E\left[T(1,1)+2p_{d}(0)v_{N-1}^{T}B^{T}S_{N}Ax_{N-1}\right.\\
 & \left.+2\bar{P}_{d}(0)u_{N-1|u_{N-1}\neq v_{N-1}}^{T}B^{T}S_{N}Ax_{N-1}+\sum_{i=0}^{N-1}x_{i}^{T}Qx_{i}\right],
\end{align*}
From this (\ref{eq:costleft2}), can be expressed as

\begin{eqnarray}
J_{N}^{*} & = & E\left[T(1,1)+2\bar{P}_{d}(0)u_{N-1|u_{N-1}\neq v_{N-1}}^{T}B^{T}S_{N}Ax_{N-1}+\sum_{i=0}^{N-1}x_{i}^{T}Qx_{i}\right.\nonumber \\
 &  & \left.+2p_{d}(0)v_{N-1}^{T}B^{T}S_{N}Ax_{N-1}+2\bar{P}_{d}(0)u_{N|u_{N}\neq v_{N}}^{T}B^{T}S_{N+1}Ax_{N}\right].\label{eq:error-11}
\end{eqnarray}
Having eliminated $x_{N}$ in (\ref{eq:TNsak}), we still have one term that depends on $x_{N}$ in (\ref{eq:costleft2}), namely
\begin{eqnarray}
E\left[\bar{P}_{d}(0)u_{N|u_{N}\neq v_{N}}^{T}B^{T}S_{N+1}Ax_{N}\right] & = & E\left[\bar{P}_{d}(0)u_{N|u_{N}\neq v_{N}}^{T}B^{T}S_{N+1}A(Ax_{N-1}+Bu_{N-1})\right]\nonumber \\
 & = & E\left[\bar{P}_{d}(0)u_{N|u_{N}\neq v_{N}}^{T}B^{T}S_{N+1}ABu_{N-1}\right.\nonumber \\
 &  & \left.+\bar{P}_{d}(0)u_{N|u_{N}\neq v_{N}}^{T}B^{T}S_{N+1}A^{2}x_{N-1}\right].\nonumber \\
 & = & E\left[\bar{P}_{d}(0)u_{N|u_{N}\neq v_{N}}^{T}B^{T}S_{N+1}ABu_{N-1}\right.\nonumber \\
 &  & +P(1)\bar{P}_{d}(0)v_{N-1}^{T}B^{T}S_{N+1}A^{2}x_{N-1}\nonumber \\
 &  & \left.+\bar{P}(1)\bar{P}_{d}(0)u_{N|u_{N}\neq v_{N},v_{N-1}}^{T}B^{T}S_{N+1}A^{2}x_{N-1}\right]\nonumber \\
 & = & E\left[\bar{P}_{d}(0)u_{N|u_{N}\neq v_{N}}^{T}B^{T}S_{N+1}ABu_{N-1}\right.\nonumber \\
 &  & +p_{d}(1)v_{N-1}^{T}B^{T}S_{N+1}A^{2}x_{N-1}\nonumber \\
 &  & \left.+\bar{P}_{d}(1)u_{N|u_{N}\neq v_{N},v_{N-1}}^{T}B^{T}S_{N+1}A^{2}x_{N-1}\right].\label{eq:sakenexpanded}
\end{eqnarray}

To evaluate $E\left[u_{N|u_{N}\neq v_{N}}^{T}B^{T}S_{N+1}ABu_{N-1}\right]$
we use a general expression which is derived in \ref{prooflong}:

\begin{align}
&E\left[u_{N|u_{N}\neq v_{N,...,N-i}}^{T}Qu_{N-k|u_{N-k}\neq v_{N-k,...,N-i}}\right] \nonumber \\
\nonumber \\
=&E\left[P(i-k+1)v_{N-i-1}^{T}Qv_{N-i-1} +v_{N-i-1}^{T}Q(P(i+1)-P(i-k+1))\right.\nonumber \\
   & \times\left(\sum_{j=1}^{N-i-\tau_{N-i-2}-2}\frac{P(i-k+j+1)-P(j-1)}{\bar{P}(j-1)}\prod_{h=2}^{j}(\frac{\bar{P}(i-k+h)}{\bar{P}(h-2)})v_{N-i-j-1}\right.\nonumber \\
   & \left.+\prod_{h=1}^{N-i-\tau_{N-i-2}-2}\frac{\bar{P}(i-k+h+1)}{\bar{P}(h-1)})v_{\tau_{N-i-2}}\right)\nonumber \\
   &\left. +\bar{P}(i+1)u_{N|u_{N}\neq v_{N},...,v_{N-i-1}}^{T}Q
    u_{N-k|u_{N-k}\neq v_{N-k},...,v_{N-i-1}}\right].\nonumber \\
    \label{eq:Uanbetingadl=0000F6st-1-1-1}
\end{align}
Note that by definition $\prod_{h}^{g}=1$ if $h>g$.
With $i=0$ , $k=1$ and  $Q=B^{T}S_{N+1}AB$ we get

\begin{align*}
&E\left[\bar{P}_{d}(0)u_{N|u_{N}\neq v_{N}}^{T}B^{T}S_{N+1}ABu_{N-1}\right]\nonumber \\
\nonumber \\
 = & E\left[\bar{P}_{d}(0)P(0)v_{N-1}B^{T}S_{N+1}ABv_{N-1}\right.
   +\bar{P}_{d}(0)p(1)v_{N-1}B^{T}S_{N+1}AB\\
   & \times\left(\sum_{j=1}^{N-\tau_{N-2}-2}\frac{p(j)}{\bar{P}(0)}v_{N-j-1}^{T}+\frac{\bar{P}(N-\tau_{N-2}-2)}{\bar{P}(0)}v_{\tau_{N-2}}^{T}\right)\\
   & \left.+\bar{P}_{d}(1)u_{N|u_{N}\neq v_{N},v_{N-1}}^{T}B^{T}S_{N+1}ABu_{N-1|u_{N-1}\neq v_{N-1}}\right].
\end{align*}
Inserting this in (\ref{eq:sakenexpanded}), we have now completely eliminated $x_{N}$ from (\ref{eq:TNsak}) and thus derived $J_{N-1}$

\begin{align}
J_{N-1}=E & \bigg[v_{N-1}^{T}T_{c}(1,1)v_{N-1}+2\bar{P}_{d}(0)P(0)v_{N-1}^{T}B^{T}S_{N+1}ABv_{N-1}\nonumber \\
& +2v_{N-1}^{T}K_{\zeta}(1,\tau_{N-2},1) \nonumber\\
 & +2v_{N-1}^{T}M(N-1)x_{N-1}+2\bar{P}_{d}(0)u_{N-1|u_{N-1}\neq v_{N-1}}^{T}B^{T}S_{N}Ax_{N-1}^{T}\nonumber \\
 & +2\bar{P}_{d}(1)u_{N|u_{N}\neq v_{N},v_{N-1}}^{T}B^{T}S_{N+1}A^{2}x_{N-1}\nonumber \\
 & +2\bar{P}_{d}(1)u_{N|u_{N}\neq v_{N},v_{N-1}}^{T}B^{T}S_{N+1}ABu_{N-1|\neq v{N-1}}\nonumber \\
 & \left.+T(1,2)+\sum_{i=0}^{N-1}x_{i}^{T}Qx_{i}\right],\label{eq:JN-1final}
\end{align}
where
\begin{align}
M(N-k)=&\sum_{i=N-k}^{N}p_{d}(i-(N-k))B^{T}S_{i+1}(A)^{i-(N-k)+1}\nonumber\\
K_{\zeta}(k,\tau,b)=& \sum_{i=N-k+1}^{N}\sum_{j=N-k+1}^{i}\bar{P}_{d}(i-(N-k+1))B^{T}S_{i+1}A^{i-j+1}\nonumber\\
 &   \times B(P(i-N+k)-P(k+j-N-1))\nonumber\\
 &   \times\left(\sum_{t=b}^{N-k-1-\tau}\frac{P(k+j-N-1+t)-P(t-1)}{\bar{P}(t-1)}\right.\nonumber\\
 &   \times\prod_{h=2}^{t}(\frac{\bar{P}(k+j-N+h-2)}{\bar{P}(h-2)})v_{N-k-t}\nonumber\\
 &   \left.+\prod_{h=1}^{N-k-1-\tau}(\frac{\bar{P}(k+j-N+h-1)}{\bar{P}(h-1)})v_{\tau}\right)\label{Kzeta}
\end{align}
and we have removed the parts that contain $v_{N-1}$ from $T(1,1)$
using

\begin{eqnarray*}
T(k,b) & = & T(k,b+1)+v_{N+1-b}^{T}T_{c}(k,b)v_{N+1-b}\\
T_{c}(k,b) & = & T_{c}(k-1,b+1)+p_{d}(b-1)B^{T}S_{N-k+1}B
\end{eqnarray*}
with $T_{c}(0,b)$ defined in (\ref{TC0}). 

The task is now to find the $v_{N-1}$ that minimizes (\ref{eq:JN-1final}).  By grouping the parts that depend on $v_{N-1}$ (\ref{eq:JN-1final})
can now be written as

\begin{align}
J_{N-1}= & E\bigg[v_{N-1}^{T}A_{11}(N-1)v_{N-1}+2v_{N-1}^{T}A_{12}(N-1,\tau_{N-2})\zeta_{N-2}\nonumber \\
 & +2\bar{P}_{d}(0)u_{N-1|u_{N-1}\neq v_{N-1}}^{T}B^{T}S_{N}Ax_{N-1}^{T}\nonumber \\
 & +2\bar{P}_{d}(1)u_{N|u_{N}\neq v_{N},v_{N-1}}^{T}B^{T}S_{N+1}A^{2}x_{N-1}+T(1,2)\nonumber \\
 & \left.+2\bar{P}_{d}(1)u_{N|u_{N}\neq v_{N},v_{N-1}}^{T}B^{T}S_{N+1}ABu_{N-1|\neq v_{N-1}}+E\sum_{i=0}^{N-1}x_{i}^{T}Qx_{i}\right]\label{eq:JN-1final-1}
\end{align}
where

\begin{eqnarray*}
A_{11}(N-1) & = & 2\bar{P}_{d}(0)P(0)B^{T}S_{N+1}AB+T_{c}(1,1)\\
A_{12}(N-1,\tau_{N-2})\zeta_{N-2} & = & M(N-1)x_{N-1}+K_{\zeta}(1,\tau_{N-2},1)\\
\zeta_{N-2} & = & \left[\begin{array}{c}
v_{N-2}\\
v_{N-3}\\
\vdots\\
v_{\tau_{N-2}}\\
x_{N-1}
\end{array}\right]
\end{eqnarray*}

This is minimized by

\[
v_{N-1}=-A_{11}^{-1}(N-1)A_{12}\zeta_{N-2}.
\]

Having inserted this into this into (\ref{eq:JN-1final-1}), what then remains to minimize is

\begin{align}
J_{N-1}^{*}= & E\bigg[-\zeta_{N-2}^{T}A_{12}^{T}(N-1)A_{11}^{-1}(N-1)A_{12}(N-1)\zeta_{N-2}\nonumber \\
 & +2\bar{P}_{d}(0)u_{N-1|u_{N-1}\neq v_{N-1}}^{T}B^{T}S_{N}Ax_{N-1}^{T}\nonumber \\
 & +2\bar{P}_{d}(1)u_{N|u_{N}\neq v_{N},v_{N-1}}^{T}B^{T}S_{N+1}A^{2}x_{N-1}+T(1,2)\nonumber \\
 & \left.+2\bar{P}_{d}(1)u_{N|u_{N}\neq v_{N},v_{N-1}}^{T}B^{T}S_{N+1}ABu_{N-1|\neq v_{N-1}}+\sum_{i=0}^{N-1}x_{i}^{T}Qx_{i}\right]\label{eq:JN-1complete-1}
\end{align}
which we now want to express in terms of variables known at time $N-2$,
i.e. we wish to find $J_{N-2}$.

We start by examining

\begin{align*}
E\left[\zeta_{N-2}^{T}A_{12}^{T}(N-1)A_{11}^{-1}(N-1)A_{12}(N-1)\zeta_{\tau_{N-2}}\right]
\end{align*}


\begin{align}
=&E\left[(K_{\eta}(1,\tau_{N-2},1)+M(N-1)x_{N-1})^{T}\right.\nonumber\\
&\left.\times A_{11}^{_{-1}}(N-1)(K_{\eta}(1,\tau_{N-2},1)+M(N-1)x_{N-1})\right]\nonumber\\\
 = &E\left[ K_{\eta}^{T}(1,\tau_{N-2},1)A_{11}^{-1}(N-1)K_{\eta}(1,\tau_{N-2},1)\nonumber \right.\\
   &\left. +x_{N-1}^{T}M^{T}(N-1)A_{11}^{-1}(N-1)M(N-1)x_{N-1}+2K_{\theta}(1,\tau_{N-2},1)x_{N-1}\right],  \nonumber \\
   \label{eq:N-2a12stuff2}
\end{align}
where

\begin{eqnarray*}
K_{\eta}(1,\tau_{},1) & = & K_{\zeta}(1,\tau_{},1)\\
K_{\theta}(1,\tau_{},1) & = & K_{\eta}^{T}(1,\tau_{},1)A_{11}^{-1}(N-1)M(N-1).
\end{eqnarray*}

Now $K_{\eta}$, $K_{\theta}$ and $K_{\zeta}$ are part of a set
of functions $\{K_{(.)}\},$ which will all have certain properties
and sub-functions. These sub-functions will be presented here for
$K_{\zeta}(k,\tau_{},b)$, but similar functions exist for all $K$
functions, which can be found in \ref{sec:conclusions}.

Now $K_{\zeta}(k,\tau_{},b)$ is defined as

\begin{eqnarray*}
K_{\zeta}(k,\tau_{},b) & = & \sum_{i=N-k+1}^{N}\sum_{j=N-k+1}^{i}\bar{P}_{d}(i-(N-k+1))B^{T}S_{i+1}A^{i-j+1}B\\
 &  & \times(P(i-N+k)-P(k+j-N-1))\\
 &  & \times\left(\sum_{t=b}^{N-k-1-\tau_{}}\frac{P(k+j-N-1+t)-P(t-1)}{\bar{P}(t-1)}\right.\\
 &  & \times\prod_{h=2}^{t}(\frac{\bar{P}(k+j-N+h-2)}{\bar{P}(h-2)})v_{N-k-t}\\
 &  & \left.+\prod_{h=1}^{N-k-1-\tau_{}}(\frac{\bar{P}(k+j-N+h-1)}{\bar{P}(h-1)})v_{\tau_{}}\right). 
\end{eqnarray*}
As can be seen, the function contains control signals from $\tau_{}$
to $N-k-b$ . Now if $\tau_{}=N-k-b$, all that remains in front of $v_{N-k-b}$ is what we name the $RL$ sub-function

\begin{align}
K_{\zeta RL}(k,b) 
  = & \sum_{i=N-k+1}^{N}\sum_{j=N-k+1}^{i}\bar{P}_{d}(i-(N-k+1))B^{T}S_{i+1}A^{i-j+1}B\nonumber\\
  & \times(P(i-N+k)-P(k+j-N-1))\prod_{h=1}^{b-1}(\frac{\bar{P}(k+j-N+h-1)}{\bar{P}(h-1)}) .\label{RLtmp}
\end{align}

If $\tau<N-k-b$ the terms affecting the last control signal
can be extracted using sub-functions with subscripts ending with a
$C$. For example,

\begin{eqnarray*}
K_{\zeta}(k,\tau_{},b) & = & K_{\zeta}(k,\tau_{},b+1)+K_{\zeta C}(k,b)v_{N-k-b}\\
K_{\zeta C}(k,b) & = & \sum_{i=N-k+1}^{N}\sum_{j=N-k+1}^{i}\bar{P}_{d}(i-(N-k+1))B^{T}S_{i+1}A^{i-j+1}\\
 &  & \times B(P(i-N+k)-P(k+j-N-1))\\
 &  & \times\frac{P(k+j-N-1+b)-P(b-1)}{\bar{P}(b-1)}\\
 &  & \times\prod_{h=2}^{b}(\frac{\bar{P}(k+j-N+h-2)}{\bar{P}(h-2)})v_{N-k-b},
\end{eqnarray*}
where $K_{\zeta}(k,\tau_{},b+1)$ now contains control signals
between $v_{N-k-b-1}$ and $v_{\tau_{}}$ but does not contain
the control signal $v_{N-k-b}$. As previously mentioned, these
sub-functions exist for all $K$ functions. For example, for $K_{\eta}$
we have
\begin{equation}
K_{\eta}(k,\tau_{},b)=K_{\eta}(k,\tau_{},b+1)+K_{\eta C}(k,b)v_{N-k-b}\label{eq:Ksaken}
\end{equation}
when $\tau_{}<N-k-b$, and since $K_{\eta}(1,\tau_{},1) = K_{\zeta}(1,\tau_{},1)$ we have
\begin{eqnarray}
K_{\eta C}(1,1) & = & K_{\zeta C}(1,1)\nonumber\\
K_{\eta RL}(1,1) & = & K_{\zeta RL}(1,1).\label{Kfuncend}
\end{eqnarray}



So returning to (\ref{eq:N-2a12stuff2}) this can be expressed in signals available at time $N-2$ using the general expression (\ref{SM2final}), with $k=1$, derived in \ref{kthings},

\begin{align}
& \left[K_{\eta}^{T}(1,\tau_{N-2},1)A_{11}^{-1}(N-1)K_{\eta}(1,\tau_{N-2},1)+\right.\nonumber\\
 & \left.x_{N-1}^{T}M^{T}(N-1)A_{11}^{-1}(N-1)M(N-1)x_{N-1}+2K_{\theta}(1,\tau_{N-2},1)x_{N-1}\right]
\nonumber\\
\nonumber\\
 =& E\left[x_{N-1}^{T}M^{T}(N-1)A_{11}^{-1}(N-1)M(N-1)x_{N-1}+2K_{ux}(2,\tau_{N-3},1)x_{N-2}\right.\nonumber \\
 & +2K_{uu}(2,\tau_{N-3},1)+2v_{N-2}^{T}K_{gx}(2,1)x_{N-2}\nonumber \\
 & +2v_{N-2}^{T}K_{gu}(2,\tau_{N-3},1,1)+K_{\alpha}(2,\tau_{N-3},1)\nonumber \\
 & \left.+v_{N-2}^{T}K_{gg}(2)v_{N-2}+2p(0)v_{N-2}^{T}K_{\theta RL}(1,1)Bv_{N-2}\right].\label{eq:afterkthings2}
\end{align}
where

\begin{eqnarray*}
K_{uu}(2,\tau_{},1) & = & \sum_{j=\tau_{}+1}^{N-3}p(N-2-j)K_{\theta}(1,j,2)Bv_{j}\\
 &  & +\bar{P}(N-3-\tau_{})K_{\theta}(1,\tau_{},2)Bv_{\tau_{}}\\
K_{ux}(2,\tau_{},1) & = & \sum_{j=\tau_{}+1}^{N-3}p(N-2-j)K_{\theta}(1,j,2)A\\
 &  & +\bar{P}(N-3-\tau_{})K_{\theta}(1,\tau_{},2)A\\
K_{gx}(2,1) & = & \bar{p}(0)K_{\theta C}(1,1)A+p(0)K_{\theta RL}(1,1)A
\end{eqnarray*}
\begin{eqnarray*}
K_{gu}(2,\tau_{},1,1) & = & \sum_{j=\tau_{}+1}^{N-3}p(N-2-j)K_{\theta C}(1,1)Bv_{j}\\
 &  & +\bar{P}(N-3-\tau_{})K_{\theta C}(1,1)KBv_{\tau_{}}\\
 &  & +K_{\beta}(2,\tau_{},1)\\
K_{\alpha}(2,\tau_{},1) & = & \sum_{j=\tau_{}+1}^{N-3}p(N-2-j)K_{\eta}^{T}(1,j,2)\\
&  & \times A_{11}^{-1}(N-1)K_{\eta}(1,j,2)\\
 &  & +\bar{P}(N-3-\tau_{})K_{\eta}^{T}(1,\tau_{},2)\\
 & & \times A_{11}^{-1}(N-1)K_{\eta}(1,\tau_{},2)\\
K_{gg}(2) & = & p(0)K_{\eta RL}^{T}(1,1)A_{11}^{-1}(N-1)K_{\eta RL}(1,1)\\
 &  & +\bar{p}(0)K_{\eta C}^{T}(1,1)A_{11}^{-1}(N-1)K_{\eta C}(1,1)\\
K_{\beta}(2,\tau_{},1) & = & \left(\sum_{j=\tau_{}+1}^{N-3}p(N-2-j)K_{\eta}(1,j,2)\right.\\
 &  & +\bar{P}(N-3-\tau_{})K_{\eta}(1,\tau_{},2)\Biggr)^{T}\\
 & & \times A_{11}^{-1}(N-1)K_{\eta C}(1,1)
\end{eqnarray*}

That finishes (\ref{eq:N-2a12stuff2}) and consequently the first term
in (\ref{eq:JN-1complete-1}) can be calculated based on signals available
at time $N-2$, except for $x_{N-1}^{T}M(N-1)A_{11}^{-1}(N-1)M^{T}(N-1)x_{N-1}$
which will be considered later. The next terms in (\ref{eq:JN-1complete-1})
we will examine are

\begin{align}
 & 2\bar{P}_{d}(0)u_{N-1|u_{N-1}\neq v_{N-1}}^{T}B^{T}S_{N}Ax_{N-1}^{T}\nonumber \\
 & +2\bar{P}_{d}(1)u_{N|u_{N}\neq v_{N},v_{N-1}}^{T}B^{T}S_{N+1}A^{2}x_{N-1}\nonumber \\
 & +2\bar{P}_{d}(1)u_{N|u_{N}\neq v_{N},v_{N-1}}^{T}B^{T}S_{N+1}ABu_{N-1|\neq v_{N-1}}\label{eq:scndpart}
\end{align}
which we will express as

\[
E\left[2F(N-1)+2\bar{P}_{d}(1)u_{N|u_{N}\neq v_{N},v_{N-1}}^{T}B^{T}S_{N+1}ABu_{N-1|\neq v_{N-1}}\right]
\]
where $F$ is defined by

\begin{align}
F(N-k)=\sum_{i=N-k}^{N}\bar{P}_{d}(i-(N-k))u_{i|u_{i}\neq v_{i},...,v_{N-k}}^{T}B^{T}S_{i+1}(A)^{i+1-(N-k)}x_{N-k}.\nonumber\\
\label{eq:FNEQ}
\end{align}

Now in (\ref{Fthings}) it is shown that

\begin{align}
E\left[F(N-k)\right]= & E\left[\sum_{i=N-k}^{N}\bar{P}_{d}(i-(N-k))u_{i|u_{i}\neq v_{i}...v_{N-k}}^{T}\right.\nonumber \\
 & \times B^{T}S_{i+1}(A)^{i+1-(N-k)}Bu_{N-k-1}\nonumber \\
 & +\sum_{i=N-k}^{N}{p}_{d}(i+1-(N-k))\nonumber \\
 & \times v_{N-k-1}^{T}B^{T}S_{i+1}(A)^{i+2-(N-k)}x_{N-k-1}\nonumber \\
 & +\sum_{i=N-k}^{N}\bar{P}_{d}(i+1-(N-k))u_{i|u_{i}\neq v_{i}...v_{N-k-1}}^{T}\nonumber \\
 & \left.\times B^{T}S_{i+1}(A)^{i+2-(N-k)}x_{N-k-1}\right].\label{eq:FNtoprove}
\end{align}

From this (\ref{eq:scndpart}) becomes:

\begin{eqnarray}
 &  & E\left[2H(2) +2\sum_{i=N-1}^{N}{p}_{d}(i+1-(N-1))v_{N-2}^{T}B^{T}S_{i+1}(A)^{i+2-(N-1)}x_{N-2}\right.\nonumber \\
 &  & \left.+2\sum_{i=N-1}^{N}\bar{P}_{d}(i+1-(N-1))u_{i|u_{i}\neq v_{i}...v_{N-2}}^{T}B^{T}S_{i+1}(A)^{i+2-(N-1)}x_{N-2}\nonumber\right]. \\
 &  & \label{eq:Fndone}
\end{eqnarray}
where the last term in (\ref{eq:scndpart}) is included in $H(2)$

\begin{align*}
H(k)= & \sum_{i=N-k+1}^{N}\bar{P}_{d}(i-(N-k+1))u_{i|u_{i}\neq v_{i},...,v_{N-k+1}}^{T}\\
 & \times B^{T}S_{i+1}\sum_{j=N-k+1}^{i}A^{i-j+1}Bu_{j-1|\neq v_{j-1},v_{N-k+1}}
\end{align*}
This expression have products of inputs we recognize from (\ref{eq:Uanbetingadl=0000F6st-1-1-1}). Using variable substitution (\ref{eq:Uanbetingadl=0000F6st-1-1-1}) can be rewritten as

\begin{align}
& E\left[u_{i|u_{i}\neq v_{i,...,N-k+1}}^{T}Qu_{j-1|u_{j-1}\neq v_{j-1,...,N-k+1}}\right]\nonumber \\
&\nonumber \\
= & E\left[P(k+j-N-1)v_{N-k}^{T}Qv_{N-k}+v_{N-k}^{T}Q\left(P(i-N+k)-P(k+j-N-1)\right)\right.\nonumber \\
 & \times\left(\sum_{t=1}^{N-k-1-\tau_{N-k-1}}\frac{P(k+j-N-1+t)-P(t-1)}{\bar{P}(t-1)}\right.\nonumber \\
 & \times \prod_{h=2}^{t}(\frac{\bar{P}(k+j-N+h-2)}{\bar{P}(h-2)})v_{N-k-t}\nonumber \\
 &\left.+\prod_{h=1}^{N-k-1-\tau_{N-k-1}}(\frac{\bar{P}(k+j-N+h-1)}{\bar{P}(h-1)})v_{\tau_{N-k-1}}\right)\nonumber \\
& \left. +\bar{P}(i-N+k)u_{i|u_{i}\neq v_{i},...,v_{N-k}}^{T}
  Qu_{j-1|u_{j-1}\neq v_{j-1},...,v_{N-k}}\right],\nonumber\\
  \label{eq:Uanbetingadl=0000F6st-1-2}
\end{align}
which allows us to derive the following expression for $H$, where
we utilise that $\tau_{N-k-1}$ is known at time $N-k$,

\begin{eqnarray*}
H(k) & = & \sum_{i=N-k+1}^{N}\sum_{j=N-k+1}^{i}\bar{P}_{d}(i-(N-k+1))\\
 &  & \times P(k+j-N-1)v_{N-k}^{T}B^{T}S_{i+1}A^{i-j+1}Bv_{N-k}\\
 &  & +\sum_{i=N-k+1}^{N}\sum_{j=N-k+1}^{i}\bar{P}_{d}(i-(N-k+1))v_{N-k}^{T}B^{T}S_{i+1}A^{i-j+1}\\
 &  & \times B(P(i-N+k)-P(k+j-N-1))\\
 &  & \times\left(\sum_{t=1}^{N-k-1-\tau_{N-k-1}}\frac{P(k+j-N-1+t)-P(t-1)}{\bar{P}(t-1)}\right.\\
 &  & \times\prod_{h=2}^{t}(\frac{\bar{P}(k+j-N+h-2)}{\bar{P}(h-2)})v_{N-k-t}\\
 &  & \left.+\prod_{h=1}^{N-k-1-\tau_{N-k-1}}(\frac{\bar{P}(k+j-N+h-1)}{\bar{P}(h-1)})v_{\tau_{N-k-1}}\right)\\
 &  & +\sum_{i=N-k+1}^{N}\sum_{j=N-k+1}^{i}\bar{P}_{d}(i-(N-k+1))\bar{P}(i-N+k)\\
 &  & \times u_{i|u_{i}\neq v_{i},...,v_{N-k}}^{T}B^{T}S_{i+1}A^{i-j+1}Bu_{j-1|u_{j-1}\neq v_{j-1},...,v_{N-k}}\\
 \end{eqnarray*}
 which we can rewrite as 
 \begin{eqnarray*}
H(k) & = & v_{N-k}^{T}K_{e}(k)v_{N-k}+v_{N-k}^{T}K_{\zeta}(k,\tau_{N-k-1},1)\\
& &+K_{a}(k),
 \end{eqnarray*}
where

\begin{eqnarray*}
K_{e}(k) & = & \sum_{i=N-k+1}^{N}\sum_{j=N-k+1}^{i}\bar{P}_{d}(i-(N-k+1))P(k+j-N-1)\\
& & \times B^{T}S_{i+1}A^{i-j+1}B\\
K_{\zeta}(k,\tau_{N-k-1},b) & = & \sum_{i=N-k+1}^{N}\sum_{j=N-k+1}^{i}\bar{P}_{d}(i-(N-k+1))B^{T}S_{i+1}A^{i-j+1}\\
 &  & \times B(P(i-N+k)-P(k+j-N-1))\\
 &  & \times\left(\sum_{t=b}^{N-k-1-\tau_{N-k-1}}\frac{P(k+j-N-1+t)-P(t-1)}{\bar{P}(t-1)}\right.\\
 &  & \times\prod_{h=2}^{t}(\frac{\bar{P}(k+j-N+h-2)}{\bar{P}(h-2)})v_{N-k-t}\\
 &  & \left.+\prod_{h=1}^{N-k-1-\tau_{N-k-1}}(\frac{\bar{P}(k+j-N+h-1)}{\bar{P}(h-1)})v_{\tau_{N-k-1}}\right)\\
K_{a}(k) & = & \sum_{i=N-k+1}^{N}\sum_{j=N-k+1}^{i}\bar{P}_{d}(i-(N-k+1))\bar{P}(i-N+k)\\
 &  & \times u_{i|u_{i}\neq v_{i},...,v_{N-k}}^{T}B^{T}S_{i+1}A^{i-j+1}u_{j-1|u_{j-1}\neq v_{j-1},...,v_{N-k}.}
\end{eqnarray*}

What remains to examine is the last parts of  (\ref{eq:JN-1complete-1}) and a term from (\ref{eq:afterkthings2}), i.e.
\begin{eqnarray}
E\left[T(1,2)-x_{N-1}^{T}M^{T}(N-1)A_{11}^{-1}(N-1)M(N-1)x_{N-1}+E\sum_{i=0}^{N-1}x_{i}^{T}Qx_{i}\right].\nonumber \\
\label{eq:doTN-1}
\end{eqnarray}

By extracting the part of $T(1,2)$ that contains $x_{N-1}$, in
the same way as $x_{N}$ was extracted from $T(0,2)$ in (\ref{TNxremove}), (\ref{eq:doTN-1})
becomes

\[
E\left[T_{noX}(1,2)+x_{N-1}^{T}S_{N-1}x_{N-1}+E\sum_{i=0}^{N-2}x_{i}^{T}Qx_{i}\right]
\]
with

\begin{eqnarray*}
S_{N-1}&=&T_{X}(1)-M^{T}(N-1)A_{11}^{-1}(N-1)M(N-1)+Q
\end{eqnarray*}
and for $k>0$, using (\ref{T_kb})
\begin{eqnarray*}
T_{noX}(k,b) & = & T_{noX}(k-1,b+1)+\sum_{i=0}^{N-k+1-b}p_{d}(N-k-i)v_{i}^{T}B^{T}S_{N-k+1}Bv_{i}\\
T_{X}(k) & = & A^{T}S_{N-k+1}A.
\end{eqnarray*}
Furthermore as shown by (\ref{eq:xSx}),

\begin{eqnarray*}
E\left[x_{N-1}S_{N-1}x_{N-1}\right] & = & E\bigg[2p_{d}(0)v_{N-2}^{T}B^{T}S_{N-1}Ax_{N-2}+x_{N-2}^{T}A^{T}S_{N-1}Ax_{N-2}\\
 &  & +2\bar{P}_{d}(0)u_{N-2|u_{N-2}\neq v_{N-2}}^{T}B^{T}S_{N-1}Ax_{N-2}\\
 &  & \left.+\sum_{i=0}^{N-2}p_{d}(N-2-i)v_{i}^{T}B^{T}S_{N-1}Bv_{i}\right]
\end{eqnarray*}
From this (\ref{eq:doTN-1}) becomes

\begin{align}
 & E\bigg[T(2,1)+2P_{d}(0)v_{N-2}^{T}B^{T}S_{N-1}Ax_{N-2}.\nonumber \\
 & \left.+2\bar{P}_{d}(0)u_{N-2|u_{N-2}\neq v_{N-2}}^{T}B^{T}S_{N-1}Ax_{N-2}+E\sum_{i=0}^{N-2}x_{i}^{T}Qx_{i}\right].\label{eq:TN2-done}
\end{align}


We can combine this with terms from (\ref{eq:Fndone}) namely,

\begin{align}
 & 2\sum_{i=N-1}^{N}\bar{P}_{d}(i+1-(N-1))u_{i|u_{i}\neq v_{i},...,v_{N-2}}^{T}B^{T}S_{i+1}(A)^{i+2-(N-1)}x_{N-2}\nonumber\\
 & +2\sum_{i=N-1}^{N}P(i+1-(N-1))\bar{P}_{d}(i-(N-1))v_{N-2}^{T}B^{T}S_{i+1}(A)^{i+2-(N-1)}x_{N-2}\label{eq:Fndone-2} 
\end{align}
to get

\[
E\left[T(2,1)+2v_{N-2}^{T}M(N-2)x_{N-2}+2F(N-2)\right]
\]
where

\begin{eqnarray*}
M(N-k) & = & \sum_{i=N-k}^{N}p_{d}(i-N-k)B^{T}S_{i+1}(A)^{i-(N-k)+1}
\end{eqnarray*}
and $F(k)$ is defined in (\ref{eq:FNtoprove}).

Finally, (\ref{eq:JN-1complete-1}) can be entirely expressed in variables
available at time $N-2$, i.e

\begin{eqnarray*}
J_{N-2} & = & E\bigg[-v_{N-2}^{T}K_{gg}(2)v_{N-2}-2p(0)v_{N-2}^{T}K_{\theta RL}(1,1)Bv_{N-2}\\
 &  & +v_{N-2}^{T}T_{c}(2,1)v_{N-2}+2v_{N-2}^{T}K_{e}(2)v_{N-2}\\
 &  & +2v_{N-2}^{T}K_{\zeta}(2,\tau_{N-3},1)+2v_{N-2}^{T}M(N-2)x_{N-2}\\
 &  & -2v_{N-2}^{T}K_{gx}(2,1)x_{N-2}-2v_{N-2}^{T}K_{gu}(2,\tau_{N-3},1,1)\\
 &  & -2K_{ux}(2,\tau_{N-3},1)x_{N-2}-2K_{uu}(2,\tau_{N-3},1)-K_{\alpha}(2,\tau_{N-3},1)\\
 &  & \left.+2F(N-2)+T(2,2)+2K_{a}(2)+\sum_{i=0}^{N-2}x_{i}^{T}Qx_{i}\right].
\end{eqnarray*}
Now, this can all be minimized with

\[
v_{N-2}=-A_{11}^{-1}(N-2)A_{12}(N-2,\tau_{N-3})\zeta_{N-3},
\]
where

\begin{eqnarray*}
A_{11}(N-2) & = & -K_{gg}(2)-2p(0)K_{\theta RL}(1,1)B\\
 &  & +T_{c}(2,1)+2K_{e}(2)
\end{eqnarray*}
and

\begin{eqnarray*}
A_{12}(N-2,\tau_{N-3})\zeta_{N-3} & = & K_{\zeta}(2,\tau_{N-3},1)-K_{gx}(2,1)x_{N-2}\\
 &  & -K_{gu}(2,\tau_{N-3},1,1)+M(N-2)x_{N-2}.
\end{eqnarray*}
This leaves to minimize

\begin{eqnarray}
J_{N-2}^{*} & = & E\bigg[-\zeta_{N-3}^{T}A_{12}^{T}(N-2)A_{11}^{-1}(N-2)A_{12}(N-2)\zeta_{N-3}\nonumber \\
 &  & -2K_{ux}(2,\tau_{N-3},1)x_{N-2}-2K_{uu}(2,\tau_{N-3},1)\nonumber \\
 &  & -K_{\alpha}(2,\tau_{N-3},1)+2F(N-2)\nonumber \\
 &  & \left.+T(2,2)+\sum_{i=0}^{N-2}x_{i}^{T}Qx_{i}+2K_{a}(2)\right].\label{JN-2sak}
\end{eqnarray}
Now we need to eliminate terms not available at time $N-3$. We start by examining

\[
2K_{uu}(2,\tau_{N-3},1)+K_{\alpha}(2,\tau_{N-3},1)
\]
and introducing the function $K_{r}$ which satisfies

\[
K_{r}(2,\tau_{N-3},1)=K_{\alpha}(2,\tau_{N-3},1)
\]
this becomes

\begin{equation}
2K_{uu}(2,\tau_{N-3},1)+K_{r}(2,\tau_{N-3},1).\label{eq:Kremaindersolve}
\end{equation}

Now, the functions $K_{uu}$ and $K_{r}$ differ a bit from most other
$K$ functions when extracting the latest control signal (see Section \ref{sec:conclusions} for details about the sub-functions):

\begin{eqnarray*}
K_{r}(k,\tau_{},b) & = & K_{r}(k,\tau_{},b+1)\\
 &  & +v_{N-k-b}^{T}K_{rCs}(k,\tau_{},b)v_{N-k-b}\\
 &  & +2K_{rCd}(N-k,\tau_{},b+1,b)v_{N-k-b}\\
K_{uu}(k,\tau_{},b) & = & K_{uu}(k,\tau_{},b+1)+v_{N-k-b}^{T}K_{uuCs}(k,b)v_{N-k-b}\\
 &  & +K_{uuCd}(k,\tau_{},b+1,b+1)v_{N-k-b}.
\end{eqnarray*}

When resolving (\ref{eq:Kremaindersolve}) we note that on time $N-3$,
$\tau_{N-3}$ is unknown. However, we can use the fact that we will at
this time know which control signal was applied at time $N-4$, i.e
$\tau_{N-4}$ to get an expression of the expected value. 

\begin{eqnarray}
 &  & E\left[2K_{uu}(2,\tau_{N-3},1)+K_{r}(2,\tau_{N-3},1)\right]\nonumber \\
 & &\nonumber \\
 & = & E\left[\sum_{j=\tau_{N-4}+1}^{N-3}p(N-3-j)(2K_{uu}(2,j,1)+K_{r}(2,j,1))\right.\nonumber \\
 &  & +\bar{P}(N-4-\tau_{N-4})(2K_{uu}(2,\tau_{N-4},1)+K_{r}(2,\tau_{N-4},1)\bigg]\nonumber \\
 & = & E\left[\sum_{j=\tau_{N-4}+1}^{N-4}p(N-3-j)(2K_{uu}(2,j,1)+K_{r}(2,j,1))\right.\nonumber \\
 &  & +p(0)v_{N-3}^{T}(2K_{uuRL}(2,1)+K_{rRL}(2,1))v_{N-3}\nonumber \\
 &  & +\bar{P}(N-4-\tau_{N-4})(2K_{uu}(2,\tau_{N-4},1)+K_{r}(2,\tau_{N-4},1)\bigg]\nonumber \\
 & = & E\left[2\sum_{j=\tau_{N-4}+1}^{N-4}p(N-3-j)K_{uu}(2,j,2))+2\bar{p}(0)v_{N-3}^{T}K_{uuCs}(2,1)v_{N-3}\right.\nonumber \\
 &  & +2\bar{P}(N-4-\tau_{N-4})K_{uu}(2,\tau_{N-4},2)\nonumber \\
 &  & +2\bar{P}(N-4-\tau_{N-4})K_{uuCd}(2,\tau_{N-4},2,2)v_{N-3}\nonumber \\
 &  & +2\sum_{j=\tau_{N-4}+1}^{N-4}p(N-3-j)K_{uuCd}(2,j,2,2)v_{N-3}\nonumber \\
 &  & +p(0)v_{N-3}^{T}\left(2K_{uuRL}(2,1)+K_{rRL}(2,1)\right)v_{N-3}\nonumber \\
 &  & +\sum_{j=\tau_{N-4}+1}^{N-4}p(N-3-j)\left(K_{r}(2,j,2)+2K_{rCd}(2,j,2,1)\right)\nonumber \\
 &  & +\bar{p}(0)v_{N-3}^{T}K_{rCs}(2,1)v_{N-3}+\bar{P}(N-4-\tau_{N-4})\bigg(K_{r}(2,\tau_{N-4},2)\nonumber \\
 &  & +2K_{rCd}(2,\tau_{N-4},2,1)v_{N-3}^{T}\bigg)\bigg]\nonumber\\
 \label{eq:Kremaindereverything}
\end{eqnarray}

By this, the terms with $K_{uu}$ and $K_{\alpha}$ can be expressed
in terms of signals available at time $N-3$. The first two terms
in (\ref{JN-2sak}) 

\begin{align*}
 &E\left[\zeta_{N-2}A_{12}^{T}(N-2,\tau_{N-3})A_{11}^{-1}(N-2)A_{12}(N-2,\tau_{N-3})\zeta_{N-2}\right.
 \\ &\left.-2K_{ux}(2,\tau_{N-3},1)x_{N-2}\right],
\end{align*}
can be expressed as

\begin{align*}
 & E\left[\left(K_{\zeta}(2,\tau_{N-3},1)-K_{gx}(2,1)x_{N-2}-K_{gu}(2,\tau_{N-3},1,1)+M(N-2)x_{N-2}\right)^{T}\right.\\
 & \times A_{11}^{-1}(N-2)\left(K_{\zeta}(2,\tau_{N-3},1)-K_{gx}(2,1)x_{N-2}-K_{gu}(2,\tau_{N-3},1,1)\right.\\
 &\left.+M(N-2)x_{N-2}\right) \left.-2K_{ux}(2,\tau_{N-3},1)x_{N-2}\right].\\
 \end{align*}


By defining for $k>1$

\begin{eqnarray*}
K_{\eta}(k,\tau_{},b) & = & K_{\zeta}(k,\tau_{},b)-K_{gu}(k,\tau_{},b,1)\\
K_{\theta}(k,\tau,b) & = & K_{\eta}(k,\tau,b)^{T}A_{11}^{-1}(N-k)\left(-K_{gx}(k,1)+M(N-k)\right)\\
 &  & +K_{ux}(k,\tau,b)\\
\end{eqnarray*}
this can be written in a similar form as before, i.e

\begin{align*}
 & E\left[K_{\eta}^{T}(2,\tau_{N-3},1)A_{11}^{-1}(N-2)K_{\eta}(2,\tau_{N-3},1)\right.\\
 & +x_{N-2}^{T}(M(N-2)-K_{gx}(2,1))^{T}A_{11}^{-1}(N-2)(M(N-2)\\
 & \left.-K_{gx}(2,1))x_{N-2}+2K_{\theta}(2,\tau_{N-3},1)x_{N-2}\right].
\end{align*}

This can be rewritten so it only depends on signals available at $N-3$
using the derivations shown in \ref{kthings}, and if it is combined with (\ref{eq:Kremaindereverything})
we get

\begin{align*}
 & E\Bigg[\zeta_{N-2}A_{12}^{T}(N-2,\tau_{N-3})A_{11}^{-1}(N-2)A_{12}(N-2,\tau_{N-3})\zeta_{N-2}+2K_{ux}(2,\tau_{N-3},1)\\
 & +2K_{uu}(2,\tau_{N-3},1)+K_{r}(2,\tau_{N-3},1)\Bigg]=
\end{align*}

\begin{align}
 & E\left[x_{N-2}^{T}(-K_{gx}(2,1)+M(N-2))A_{11}^{-1}(-K_{gx}(2,1)+M(N-2))x_{N-2}\right.\nonumber \\
 & +2K_{ux}(3,\tau_{N-4},1)x_{N-3}+2K_{uu}(3,\tau_{N-4},1)+2v_{N-3}^{T}K_{gx}(3,1)x_{N-3}\nonumber \\
 & +2v_{N-3}^{T}K_{gu}(3,\tau_{N-4},1,1)+2p(0)v_{N-3}^{T}K_{\theta RL}(2,1)Bv_{N-3}\nonumber \\
 & \left.+v_{N-3}^{T}K_{gg}(3)v_{N-3}+K_{r}(3,\tau_{N-4},1)\right], \nonumber \\
 \label{eq:Msak}
\end{align}
with for $k>2$
\begin{align*}
K_{\alpha}(k,\tau,b)= & \sum_{j=\tau+1}^{N-k-b}p(N-k-j) K_{\eta}^{T}(k-1,j,b+1)\\
 & \times A_{11}^{-1}(N-k+1)K_{\eta}(k-1,j,b+1)\\
 & +\bar{P}(N-k-1-\tau)K_{\eta}^{T}(k-1,\tau,b+1)\\
 & \times A_{11}^{-1}(N-k+1)K_{\eta}(k-1,\tau,b+1)\\
 K_{gg}(k)= & p(0)K_{\eta RL}^{T}(k-1,1)A_{11}^{-1}(N-k+1)K_{\eta RL}(k-1,1)\\
 & +\bar{p}(0)K_{\eta C}^{T}(k-1,1) A_{11}^{-1}(N-k+1)K_{\eta C}(k-1,1)\\
 & +2\bar{p}(0)K_{uuCs}(k-1,1) +\bar{p}(0)K_{rCs}(k-1,1)\\
 & +p(0)\left(2K_{uuRL}(k-1,1)+K_{rRL}(k-1,1)\right)\\
K_{\beta}(k,\tau,b)= &\left( \sum_{j=\tau+1}^{N-k-b}p(N-k-j)K_{\eta}(k-1,j,b+1)\right.\\
 & +\bar{P}(N-k-1-\tau_{N-3})K_{\eta}(k-1,\tau,b+1)\bigg)^{T}\\
 & \times A_{11}^{-1}(N-k+1)K_{\eta C}(k-1,1)\\
 K_{uu}(k,\tau,b)= & \sum_{j=\tau+1}^{N-k-b}p(N-k-j)K_{\theta}(k-1,j,b+1)Bv_{j}\\
 & +\bar{P}(N-k-1-\tau)K_{\theta}(k-1,\tau,b+1)Bv_{\tau}\\
 & +\sum_{j=\tau+1}^{N-k-b}p(N-k-j)K_{uu}(k-1,j,b+1)\\
 & +\bar{P}(N-k-1-\tau)K_{uu}(k-1,\tau,b+1)\\
   K_{gu}(k,\tau,b,h)=& \sum_{j=\tau+1}^{N-k-b}p(N-k-j)K_{\theta C}(k-1,h)Bv_{j}\\
 & +\bar{P}(N-k-1-\tau)K_{\theta C}(k-1,h)Bv_{\tau}\\
 & +K_{\beta}(k,\tau,b)\\
 & +\sum_{j=\tau+1}^{N-k-b}p(N-k-j)K_{uuCd}(k-1,j,b+1,h+1)\\
 & +\bar{P}(N-k-1-\tau)K_{uuCd}(k-1,\tau,b+1,h+1)\\
 & +\sum_{j=\tau+1}^{N-k-b}p(N-k-j)K_{rCd}(k-1,j,b+1)\\
 & +\bar{P}(N-k-1-\tau)K_{rCd}(k-1,\tau,b+1,h)\\
\end{align*}
\begin{align*}
  K_{ux}(k,\tau,b)= & \sum_{j=\tau+1}^{N-k-b}p(N-k-j)K_{\theta}(k-1,j,b+1)A\\
 & +\bar{P}(N-k-1-\tau_{N-3})K_{\theta}(k-1,\tau,b+1)A\\
 K_{gx}(k,b)= & \bar{p}(0)K_{\theta C}(k-1,b)A+p(0)K_{\theta RL}(k-1,b)A\\
  K_{r}(k,\tau,b)=& \sum_{j=\tau+1}^{N-k-b}p(N-k-j)K_{r}(k-1,j,b+1)\\
 & +\bar{P}(N-k-1-\tau)K_{r}(k-1,\tau,b+1)\\
 & +K_{\alpha}(k,\tau,b)
\end{align*}

Now, for

\begin{align*}
 & E\bigg[-x_{N-2}^{T}(M(N-2)-K_{gx}(2,1))^{T}A_{11}^{-1}(M(N-2)-K_{gx}(2,1))x_{N-2}\\
 & \left. T(2,2)+\sum_{i=0}^{N-2}x_{i}^{T}Qx_{i}\right],
\end{align*}
from (\ref{JN-2sak}) and (\ref{eq:Msak}). Doing the same calculations
as to get (\ref{eq:TN2-done}) this becomes

\begin{align}
 & E\left[+2P_{d}(0)v_{N-3}^{T}B^{T}S_{N-2}Ax_{N-3}+2\bar{P}_{d}(0)u_{N-3|u_{N-3}\neq v_{N-3}}^{T}B^{T}S_{N-2}Ax_{N-3}\right.\nonumber \\
 & \left.T(3,1)+E\sum_{i=0}^{N-3}x_{i}^{T}Qx_{i}\right].\nonumber \\
 \label{eq:TN3-done}
\end{align}
where for $k>1$

\begin{eqnarray*}
S_{N-k} & = & -(M(N-k)-K_{gx}(k,\tau_{N-k-1}))^{T}A_{11}^{-1}(N-k)(M(N-k)-K_{gx}(k,\tau_{N-k-1}))\\
 &  & +T_{X}(k)+Q\\
\end{eqnarray*}

Now adding  $2F(N-2)+2K_{a}(2),$  from (\ref{JN-2sak}) to (\ref{eq:TN3-done})
and using the same calculations as those previously done on $H(2)$, we get

\begin{align*}
 & E\bigg[T(3,1)+2v_{N-3}^{T}M(N-3)x_{N-3}+2v_{N-3}^{T}K_{e}(3)v_{N-3}\\
 & \left.+2K_{\zeta}(3,\tau_{N-4},1)v_{N-3}+2K_{a}(3)+2F(N-3)+E\sum_{i=0}^{N-3}x_{i}^{T}Qx_{i}\right].
\end{align*}

Now we have rewritten all the terms in (\ref{JN-2sak}) so they only
depend on signals available at time $N-3$, so if we put them together
we get:

\begin{eqnarray*}
J_{N-3} & = & E\bigg[-v_{N-3}^{T}K_{gg}(3)v_{N-3}\\
 &  & -2p(0)v_{N-3}^{T}K_{\theta RL}(2,1)Bv_{N-3}\\
 &  & +v_{N-3}^{T}T_{c}(3,1)v_{N-3}+2v_{N-3}^{T}K_{e}(3)v_{N-3}\\
 &  & +2K_{\zeta}(3,\tau_{N-4},1)v_{N-3}-2v_{N-3}^{T}K_{gx}(3,1)x_{N-3}\\
 &  & -2K_{gu}(3,\tau_{N-4},1,1)v_{N-3}+2v_{N-3}^{T}M(N-3)x_{N-3}\\
 &  & -2K_{ux}(3,\tau_{N-4},1,1)x_{N-3}-2K_{uu}(3,\tau_{N-4},1)\\
 &  & -K_{r}(3,\tau_{N-4},1)+2F(N-3)\\
 &  & \left.+T(3,2)+2K_{a}(3)+E\sum_{i=0}^{N-3}x_{i}^{T}Qx_{i}\right].
\end{eqnarray*}

From this the control signal for time $N-3$ can be derived, i.e.

\[
v_{N-3}=A_{11}^{-1}(N-3)A_{12}(N-3,\tau_{N-4})\zeta_{N-k-1},
\]
where

\begin{eqnarray*}
A_{11}(N-3) & = & -K_{gg}(3)-2p(0)K_{\theta RL}(2,1)B+T_{c}(N-3)\\
& & +2K_{e}(3)\\
A_{12}(N-3,\tau_{N-4})\zeta_{N-4} & = & K_{\zeta}(3,\tau_{N-4},1)-K_{gx}(3,1)x_{N-3}\\
 &  & -K_{gu}(3,\tau_{N-k-1},1,1)+M(N-k)x_{N-3}
\end{eqnarray*}

By repeating the previous calculations an expression for $A_{11}(N-k)$
and $A_{12}(N-k)$ can be found for $k>2$, i.e.

\begin{eqnarray*}
A_{11}(N-k) & = & -K_{gg}(k,1)-2p(0)K_{\theta RL}(k-1,1)B\\
 &  & +T_{c}(k)+2K_{e}(k)\\
A_{12}(N-k,\tau_{N-k-1})\zeta_{N-k-1} & = & K_{\zeta}(k,\tau_{N-k-1},1)-K_{gx}(k,1)x_{N-k}\\
 &  & -K_{gu}(k,\tau_{N-k-1},1,1)+M(N-k)x_{N-k}
\end{eqnarray*}
and a consequent optimal feedback control

\[
v_{N-k}=A_{11}^{-1}(N-k)A_{12}(N-k,\tau_{N-k-1})\zeta_{N-k-1}.
\]

\subsection{Proof for (\ref{eq:Uanbetingadl=0000F6st-1-1-1})}\label{prooflong}

We wish to find an expression for the expected value of \\$u_{N|u_{N}\neq v_{N,...,N-i}}^{T}Qu_{N-k|u_{N-k}\neq v_{N-k,...,N-i}}$
in the case where the actuator signal $u_{N-i-2}$ is known. We start by examining the probability that the control signal applied at time $N-k$ is $v_{N-i-1}$. As we know no later signal  has arrived, the probability of this signal being applied is simply the probability it has arrived by this time, i.e $P(i-k+1)$. If this signal is applied at time $N-k$ it will also be applied at time $N$  since no later signal has reached the actuator at this time. So

\[
E\left[u_{N|u_{N}\neq v_{N},...,v_{N-i}}^{T}Qu_{N-k|u_{N-k}\neq v_{N-k},...,v_{N-i}}\right]=
\]

\begin{align}
 & E\left[P(i-k+1)v_{N-i-1}^{T}Qv_{N-i-1}\right. \nonumber \\
 &\left.+\bar{P}(i-k+1)u_{N|u_{N}\neq v_{N},...,v_{N-i}}^{T}Qu_{N-k|u_{N-k}\neq v_{N-k},...,v_{N-i-1}}\right]  \nonumber\\ \label{uanfirst-1sak}
\end{align}

Next we examine the probability that the control signal applied at time $N$ is $v_{N-i-1}$ given that this signal was not been applied at time $N-k$, i.e

\[
p(u_{N}=v_{N-i-1}|u_{N}\neq v_{N},...,v_{N-i}|u_{N-k}\neq v_{N-k},...,v_{N-i-1})
\]

Now this can be seen to be the probability that $v_{N-i-1}$ has arrived at time $N$ given that it has not yet arrived at time $N-k$, i.e the probability that the delay of $v_{N-i-1}$ is less or equal to $i+1$  given that it is greater then $i+1-k$ . This probability can be calculated as
\begin{align}
& \frac{P(i+1)-P(i-k+1)}{\bar{P}(i-k+1)}.\label{probapplied}
\end{align}



Using this, the second term of (\ref{uanfirst-1sak}) can be written as
\[
\bar{P}(i-k+1)u_{N|u_{N}\neq v_{N},...,v_{N-i}}^{T}Qu_{N-k|u_{N-k}\neq v_{N-k},...,v_{N-i-1}}=
\]

\begin{align}
 & E\left[(P(i+1)-P(i-k+1))v_{N-i-1}^{T}Qu_{N-k|u_{N-k}\neq v_{N-k},...,v_{N-i-1}}\right.\nonumber\\ 
 &\left.+\bar{P}(i+1)u_{N|u_{N}\neq v_{N},...,v_{N-i-1}}^{T}Qu_{N-k|u_{N-k}\neq v_{N-k},...,v_{N-i-1}}\right]  \label{uanfirst-1sak2}
\end{align}

Now we will evaluate, $E\left[u_{N-k|u_{N-k}\neq v_{N-k},...,v_{N-i-1}}\right]$
using the fact that we know the actuator signal at time $N-i-2$,
and consequently know which possible signals may be applied at time $N-i-1$. For a signal to be applied at time $N-i-1$ the signal must have arrived, this probability can be calculated the same way as was done in (\ref{probapplied}). Furthermore no later signal must have arrived by this time. As the delays between signals is independent this can be taken into account by taking the complimentary probability of each later control signal arriving. Combining these two requirements yields,   

\begin{align}
&E\left[u_{N-k|u_{N-k}\neq v_{N-k},...,v_{N-i-1},u_{N-i-2}=v_{\tau_{N-i-2}}}\right] 
\nonumber \\
\nonumber \\
=&E\left[\sum_{j=1}^{N-i-\tau_{N-i-2}-2}\frac{P(i-k+j+1)-P(j-1)}{\bar{P}(j-1)}\right.\nonumber \\
   & \times\prod_{h=2}^{j}(1-\frac{P(i-k+h)-P(h-2)}{\bar{P}(h-2)})v_{N-i-j-1}\nonumber \\
   & \left.+\prod_{h=1}^{N-i-\tau_{N-i-2}-2}(1-\frac{P(i-k+h+1)-P(h-1)}{\bar{P}(h-1)})v_{\tau_{N-i-2}}\right]\nonumber \\
  = & \sum_{j=1}^{N-i-\tau_{N-i-2}-2}\frac{P(i-k+j+1)-P(j-1)}{\bar{P}(j-1)}\prod_{h=2}^{j}(\frac{\bar{P}(i-k+h)}{\bar{P}(h-2)})v_{N-i-j-1}\nonumber \\
   & +\prod_{h=1}^{N-i-\tau_{N-i-2}-2}\frac{\bar{P}(i-k+h+1)}{\bar{P}(h-1)})v_{\tau_{N-i-2}},\label{eq:uansak-1thing}
\end{align}
where $\tau_{N-i-2}$ is the sample of the control signal that is applied
at time $N-i-2$.

Finally, combining (\ref{uanfirst-1sak}) , (\ref{uanfirst-1sak2}), and
(\ref{eq:uansak-1thing}) 

\begin{align}
& E\left[u_{N|u_{N}\neq v_{N},...,v_{N-i}}^{T}Qu_{N-k|u_{N-k}\neq v_{N-k},...,v_{N-i}}\right] \nonumber\\ 
=& E\left[P(i-k+1)v_{N-i-1}^{T}Qv_{N-i-1} +(P(i+1)-P(i-k+1))v_{N-i-1}^{T}Q\right. \nonumber \\
  & \times\left(\sum_{j=1}^{N-i-\tau_{N-i-2}-2}\frac{P(i-k+j+1)-P(j-1)}{\bar{P}(j-1)}\prod_{h=2}^{j}(\frac{\bar{P}(i-k+h)}{\bar{P}(h-2)})v_{N-i-j-1} \right. \nonumber \\
   & \left.+\prod_{h=1}^{N-i-\tau_{N-i-2}-2}\frac{\bar{P}(i-k+h+1)}{\bar{P}(h-1)})v_{\tau_{N-i-2}}\right)\nonumber \\
   & \left. +\bar{P}(i+1)u_{N|u_{N}\neq v_{N},...,v_{N-i-1}}^{T}Q u_{N-k|u_{N-k}\neq v_{N-k},...,v_{N-i-1}}\right]. \nonumber \\
\end{align}

\subsection{General proof used for (\ref{eq:afterkthings2})}\label{kthings}

We wish to find an expression for

\begin{eqnarray}
 &  & E\left[K_{\eta}^{T}(k,\tau_{N-k-1},1)A_{11}^{-1}(N-k)K_{\eta}(k,\tau_{N-k-1},1)\right.\nonumber \\
 &  & \left.+2K_{\theta}(k,\tau_{N-k-1},1)x_{N-k}\right]\label{eq:N-k-1a12stuff-1}
\end{eqnarray}
expressed in signals available at time $N-k-1$ and the control signal $v_{N-k-1}$.
To start, we examine $E(K_{\eta}(k,\tau_{N-k-1},1))$. At time $N-k-1$,
$\tau_{N-k-1}$ is unknown. However, we can use the fact that we will at
this time know which control signal was applied at time $N-k-2$, i.e
$\tau_{N-k-2}$, to get an expression of the expected value. This is a similar problem to that of (\ref{eq:uansak-1thing}), and by applying the same reasoning as was done there, this expected value can be derived as

\begin{align}
&E(K_{\eta}(k,\tau_{N-k-1},1)) \nonumber\\ 
=& \sum_{j=\tau_{N-k-2}+1}^{N-k-1}\frac{p(N-k-1-j)}{\bar{P}(N-k-2-j)}\prod_{h=1}^{N-k-1-j}(1-\frac{p(h-1)}{\bar{P}(h-2)})K_{\eta}(k,j,1)\nonumber\\ 
   & +\prod_{h=1}^{N-k-1-\tau_{N-k-2}}(1-\frac{p(h-1)}{\bar{P}(h-2)})K_{\eta}(k,\tau_{N-k-2},1)\nonumber\\ 
  = & \sum_{j=\tau_{N-k-2}+1}^{N-k-1}\frac{p(N-k-1-j)}{\bar{P}(N-k-2-j)}\nonumber\\ 
   & \times\prod_{h=1}^{N-k-1-j}(\frac{\bar{P}(h-2)}{\bar{P}(h-2)}-\frac{p(h-1)}{\bar{P}(h-2)})K_{\eta}(k,j,1)\nonumber\\ 
   & +\prod_{h=1}^{N-k-1-\tau_{N-k-2}}(\frac{\bar{P}(h-2)}{\bar{P}(h-2)}-\frac{p(h-1)}{\bar{P}(h-2)})K_{\eta}(k,\tau_{N-k-2},1)\nonumber\\ 
  = & \sum_{j=\tau_{N-k-2}+1}^{N-k-1}\frac{p(N-k-1-j)}{\bar{P}(N-k-2-j)}\prod_{h=1}^{N-k-1-j}(\frac{\bar{P}(h-1)}{\bar{P}(h-2)})K_{\eta}(k,j,1)\nonumber\\ 
   & +\prod_{h=1}^{N-k-1-\tau_{N-k-2}}(\frac{\bar{P}(h-1)}{\bar{P}(h-2)})K_{\eta}(k,\tau_{N-k-2},1)\nonumber\\ 
  = & \sum_{j=\tau_{N-k-2}+1}^{N-k-1}\frac{p(N-k-1-j)}{\bar{P}(N-k-2-j)}\bar{P}(N-k-2-j)K_{\eta}(k,j,1)\nonumber\\ 
   & +\bar{P}(N-k-2-\tau_{N-k-2})K_{\eta}(k,\tau_{N-k-2},1)\nonumber\\ 
  = & \sum_{j=\tau_{N-k-2}+1}^{N-k-1}p(N-k-1-j)K_{\eta}(k,j,1)\nonumber\\ 
   & +\bar{P}(N-k-2-\tau_{N-k-2})K_{\eta}(k,\tau_{N-k-2},1).\label{Simpthing}
\end{align}

From this and using the $K$ subfunctions as expressed in (\ref{eq:Ksaken}) and the definition of RL-functions as defined above (\ref{RLtmp}), the first term in (\ref{eq:N-k-1a12stuff-1})
becomes

\begin{align*}
 & E\left[\sum_{j=\tau_{N-k-2}+1}^{N-k-2}p(N-k-1-j)\left((K_{\eta}(k,j,2)+K_{\eta C}(k,1)v_{N-k-1})^{T}A_{11}^{-1}(N-k)\right.\right.\\
 & \left.\times(K_{\eta}(k,j,2)+K_{\eta C}(k,1)v_{N-k-1})\right)\\
 & +p(0)v_{N-k-1}^{T}K_{\eta RL}^{T}(k,1)A_{11}^{-1}(N-k)K_{\eta RL}(k,1)v_{N-k-1}\\
 & +\bar{P}(N-k-2-\tau_{N-k-2})(K_{\eta}(k,\tau_{N-k-2},2)+K_{\eta C}(k,1)v_{N-k-1})^{T}A_{11}^{-1}(N-k)\\
 & \times(K_{\eta}(k,\tau_{N-k-2},2)+K_{\eta C}(k,1)v_{N-k-1})\bigg],
\end{align*}

\begin{align}
= & E\left[\sum_{j=\tau_{N-k-2}+1}^{N-k-2}p(N-k-1-j)K_{\eta}^{T}(k,j,2)A_{11}^{-1}(N-k)K_{\eta}(k,j,2)\right.\nonumber \\
 & +\bar{P}(N-k-2-\tau_{N-k-2})K_{\eta}^{T}(k,\tau_{N-k-2},2)A_{11}^{-1}(N-k)K_{\eta}(k,\tau_{N-k-2},2)\nonumber \\
 & +p(0)v_{N-k-1}^{T}K_{\eta RL}^{T}(k,1)A_{11}^{-1}(N-k)K_{\eta RL}(k,1)v_{N-k-1}\nonumber \\
 &+\bar{p}(0) v_{N-k-1}^{T}K_{\eta C}^{T}(k,1)A_{11}^{-1}(N-k)K_{\eta C}(k,1)v_{N-k-1}\nonumber \\
 & +2\left(\sum_{j=\tau_{N-k-2}+1}^{N-k-2}p(N-k-1-j)K_{\eta}(k,j,2)\right.\nonumber \\
 & +\bar{P}(N-k-2-\tau_{N-k-2})K_{\eta}(k,\tau_{N-k-2},2)\bigg)^{T}
 A_{11}^{-1}(N-k)K_{\eta C}(k,1)v_{N-k-1}\bigg].\nonumber \\
 \label{eq:Ksaveusesak-1}
\end{align}

Now, for $E\left[2K_{\theta}(k,\tau_{N-k-1},1)x_{N-k}\right]$ in (\ref{eq:N-k-1a12stuff-1})
, we have,

	\begin{eqnarray}	
	K_{\theta}(k,\tau_{N-k-1},1)x_{N-k} & = & K_{\theta}(k,\tau_{N-k-1},1)Ax_{N-k-1}\nonumber \\	
	&&+K_{\theta}(k,\tau_{N-k-1},1)Bu_{N-k-1}.\label{eq:diffxall-1}	
	\end{eqnarray}
	
For $E(K_{\theta}(k,\tau_{N-k-1},1)Bu_{N-k-1}) $ we will first do the same calculations as in (\ref{Simpthing}) using the fact that we know which control signals were applied at time $N-k-2$, which yields

\begin{align}
& E(K_{\theta}(k,\tau_{N-k-1},1)Bu_{N-k-1}) \nonumber \\
\nonumber \\
 = & E\left[\sum_{j=\tau_{N-k-2}+1}^{N-k-1}p(N-k-1-j)K_{\theta}(k,j,1)Bv_{j}\right.\nonumber \\
   & +\bar{P}(N-k-2-\tau_{N-k-2})K_{\theta}(k,\tau_{N-k-2},1)Bv_{\tau_{N-k-2}}\bigg]\nonumber \\
 = & E\left[\sum_{j=\tau_{N-k-2}+1}^{N-k-2}p(N-k-1-j)K_{\theta}(k,j,1)Bv_{j}\right.\nonumber \\
  & +\bar{P}(N-k-2-\tau_{N-k-2})K_{\theta}(k,\tau_{N-k-2},1)Bv_{\tau_{N-k-2}}\nonumber \\
   & +p(0)v_{N-k-1}^{T}K_{\theta RL}(k,1)Bv_{N-k-1}\bigg]\nonumber \\
  = & E\left[\sum_{j=\tau_{N-k-2}+1}^{N-k-2}p(N-k-1-j)K_{\theta}(k,j,2)Bv_{j}\right.\nonumber \\
   & +\bar{P}(N-k-2-\tau_{N-k-2})K_{\theta}(k,\tau_{N-k-2},2)Bv_{\tau_{N-k-2}}\nonumber \\
   & +\bar{P}(N-k-2-\tau_{N-k-2})v_{N-k-1}^{T}K_{\theta C}(k,1)Bv_{\tau_{N-k-2}})\nonumber \\
   & +\sum_{j=\tau_{N-k-2}+1}^{N-k-2}p(N-k-1-j)v_{N-k-1}^{T}K_{\theta C}(k,1)Bv_{j}\nonumber \\
   & +p(0)v_{N-k-1}^{T}K_{\theta RL}(k,1)Bv_{N-k-1}\bigg]\label{diffxu-1}
\end{align}

Now if we do the same for $E\left[K_{\theta}(k,\tau_{N-k-1},1)Ax_{N-k-1}\right] $ and use that $\bar{P}(N-k-2-\tau_{N-k-2})=1-\sum_{j=0}^{N-k-2-\tau_{N-k-2}}p(j)$, it becomes,

\begin{align}
& E\left[K_{\theta}(k,\tau_{N-k-1},1)Ax_{N-k-1}\right] \nonumber \\
\nonumber \\
= & E\left[\sum_{j=\tau_{N-k-2}+1}^{N-k-1}p(N-k-1-j)K_{\theta}(k,j,1)Ax_{N-k-1}\right.\nonumber \\
  & +\bar{P}(N-k-2-\tau_{N-k-2})K_{\theta}(k,\tau_{N-k-2},1)Ax_{N-k-1}\bigg]\nonumber \\
 = & E\left[\sum_{j=\tau_{N-k-2}+1}^{N-k-2}p(N-k-1-j)K_{\theta}(k,j,2)Ax_{N-k-1}\right.\nonumber \\
 &+p(0)v_{N-k-1}^{T}K_{\theta RL}(k,1)Ax_{N-k-1}
   \nonumber \\
   & +\bar{p}(0)v_{N-k-1}^{T}K_{\theta C}(k,1)Ax_{N-k-1}\nonumber \\
   & +\bar{P}(N-k-2-\tau_{N-k-2})K_{\theta}(k,\tau_{N-k-2},2)Ax_{N-k-1} \bigg].\label{eq:diffxx-1}
\end{align}



So, from (\ref{eq:Ksaveusesak-1}), (\ref{diffxu-1}) and (\ref{eq:diffxx-1}) , (\ref{eq:N-k-1a12stuff-1}) can be expressed in signals available at time $N-k-1$.

\begin{align}
= & E\left[\sum_{j=\tau_{N-k-2}+1}^{N-k-2}p(N-k-1-j)K_{\eta}^{T}(k,j,2)A_{11}^{-1}K_{\eta}(k,j,2)\right.\nonumber \\
 & +\bar{P}(N-k-2-\tau_{N-k-2})K_{\eta}^{T}(k,\tau_{N-k-2},2)A_{11}^{-1}K_{\eta}(k,\tau_{N-k-2},2)\nonumber \\
 & +p(0)v_{N-k-1}^{T}K_{\eta RL}^{T}(k,1)A_{11}^{-1}K_{\eta RL}(k,1)v_{N-k-1}\nonumber \\
 &+\bar{p}(0) v_{N-k-1}^{T}K_{\eta C}^{T}(k,1)A_{11}^{-1}K_{\eta C}(k,1)v_{N-k-1}\nonumber \\
 & +2\left(\sum_{j=\tau_{N-k-2}+1}^{N-k-2}p(N-k-1-j)K_{\eta}(k,j,2)\right.\nonumber \\
 &\left. +\bar{P}(N-k-2-\tau_{N-k-2})K_{\eta}(k,\tau_{N-k-2},2)\right)^{T}
 A_{11}^{-1}K_{\eta C}(k,1)v_{N-k-1}\nonumber \\
   & +2\sum_{j=\tau_{N-k-2}+1}^{N-k-2}p(N-k-1-j)K_{\theta}(k,j,2)Bv_{j}\nonumber \\
   & +2\bar{P}(N-k-2-\tau_{N-k-2})K_{\theta}(k,\tau_{N-k-2},2)Bv_{\tau_{N-k-2}}\nonumber \\
   & +2\bar{P}(N-k-2-\tau_{N-k-2})v_{N-k-1}^{T}K_{\theta C}(k,1)Bv_{\tau_{N-k-2}})\nonumber \\
   & +2\sum_{j=\tau_{N-k-2}+1}^{N-k-2}p(N-k-1-j)v_{N-k-1}^{T}K_{\theta C}(k,1)Bv_{j}\nonumber \\
   & +2p(0)v_{N-k-1}^{T}K_{\theta RL}(k,1)Bv_{N-k-1}\nonumber \\
   & +2p(0)v_{N-k-1}^{T}K_{\theta RL}(k,1)Ax_{N-k-1}\nonumber \\
   & +2\sum_{j=\tau_{N-k-2}+1}^{N-k-2}p(N-k-1-j)K_{\theta}(k,j,2)Ax_{N-k-1}\nonumber \\
   & +2\bar{p}(0)v_{N-k-1}^{T}K_{\theta C}(k,1)Ax_{N-k-1}\nonumber \\
   & +2\bar{P}(N-k-2-\tau_{N-k-2})K_{\theta}(k,\tau_{N-k-2},2)Ax_{N-k-1} \bigg]\label{SM2final}
\end{align}

\subsection{Proof for (\ref{eq:FNtoprove})}
\label{Fthings}

The definition of $F$ is

\[
F(N-k)=\sum_{i=N-k}^{N}\bar{P}_{d}(i-(N-k))u_{i|u_{i}\neq v_{i}...v_{N-k}}^{T}B^{T}S_{i+1}(A)^{i+1-(N-k)}x_{N-k}
\]
and we wish to derive an expression for $E\left[F(N-k)\right]$, not dependant on the state $x_{N-k}$

\begin{eqnarray}
E\left[F(N-k)\right] & = & E\left[\sum_{i=N-k}^{N}\bar{P}_{d}(i-(N-k))u_{i|u_{i}\neq v_{i}...v_{N-k}}^{T}B^{T}S_{i+1} \right. \nonumber \\
& &\left. \times (A)^{i+1-(N-k)}(Ax_{N-k-1}+Bu_{N-k-1})\right]\nonumber \\
 & = & E\left[\sum_{i=N-k}^{N}\bar{P}_{d}(i-(N-k))u_{i|u_{i}\neq v_{i}...v_{N-k}}^{T}B^{T}S_{i+1}(A)^{i+1-(N-k)}Bu_{N-k-1}\right.\nonumber \\
 &  & \left.+\sum_{i=N-k}^{N}\bar{P}_{d}(i-(N-k))u_{i|u_{i}\neq v_{i}...v_{N-k}}^{T}B^{T}S_{i+1}(A)^{i+2-(N-k)}x_{N-k-1}\right]\nonumber \\
 \label{FNdel}
\end{eqnarray}
where

\begin{equation*}
E\left[\sum_{i=N-k}^{N}\bar{P}_{d}(i-(N-k))u_{i|u_{i}\neq v_{i}...v_{N-k}}^{T}B^{T}S_{i+1}(A)^{i+2-(N-k)}x_{N-k-1}\right]
\end{equation*}
\begin{align}
 =& E\left[\sum_{i=N-k}^{N}P(i+1-(N-k))\bar{P}_{d}(i-(N-k))v_{N-k-1}^{T}B^{T}S_{i+1}(A)^{i+2-(N-k)}x_{N-k-1}\right.\nonumber\\
 & +\sum_{i=N-k}^{N}\bar{P}(i+1-(N-k))\bar{P}_{d}(i-(N-k))u_{i|u_{i}\neq v_{i}...v_{N-k-1}}^{T} \nonumber\\&\left.\times B^{T}S_{i+1}(A)^{i+2-(N-k)}x_{N-k-1}\right].\nonumber\\\label{eq:FsakPDs}
\end{align}

Furthermore,

\begin{eqnarray*}
\bar{P}(i+1-(N-k))\bar{P}_{d}(i-(N-k)) & = & \bar{P}_{d}(i-(N-k))-P(i+1-(N-k))\bar{P}_{d}(i-(N-k))\\
 & = & \bar{P}_{d}(i-(N-k))-p_{d}(i+1-(N-k))\\
 & = & \bar{P}_{d}(i+1-(N-k))
\end{eqnarray*}
using this and (\ref{pddef3}),  (\ref{eq:FsakPDs}) simplifies to

\begin{align*}
 & E\left[\sum_{i=N-k}^{N}{p}_{d}(i+1-(N-k))v_{N-k-1}^{T}B^{T}S_{i+1}(A)^{i+2-(N-k)}x_{N-k-1}\right.\\
 & \left.+\sum_{i=N-k}^{N}\bar{P}_{d}(i+1-(N-k))u_{i|u_{i}\neq v_{i}...v_{N-k-1}}^{T}B^{T}S_{i+1}(A)^{i+2-(N-k)}x_{N-k-1}\right].\\
\end{align*}

So, 

\begin{eqnarray*}
E\left[F(N-k)\right] & = & E\left[\sum_{i=N-k}^{N}\bar{P}_{d}(i-(N-k))u_{i|u_{i}\neq v_{i}...v_{N-k}}^{T}B^{T}S_{i+1}(A)^{i+1-(N-k)}Bu_{N-k-1}\right.\\
 &  & +\sum_{i=N-k}^{N}{p}_{d}(i+1-(N-k))v_{N-k-1}^{T}B^{T}S_{i+1}(A)^{i+2-(N-k)}x_{N-k-1}\\
 &  & \left.+\sum_{i=N-k}^{N}\bar{P}_{d}(i+1-(N-k))u_{i|u_{i}\neq v_{i}...v_{N-k-1}}^{T}B^{T}S_{i+1}(A)^{i+2-(N-k)}x_{N-k-1}\right].
\end{eqnarray*}

\section{Implementation summary} 
\label{sec:conclusions}

The formulas needed to calculate the optimal control are here summarized
to facilitate an easy implementation.

Note that the only information not known prior to implementation are the previous delays $\tau_{(\cdot)}$. Thus all calculations not including $\tau_{(\cdot)}$ can be done offline or before hand.

\begin{eqnarray*}
v_{N-k} & = & A_{11}^{-1}(N-k)A_{12}(N-k)\zeta_{N-k-1}.\\
A_{11}(N) & = & T_{c}(0,1)\\
A_{12}(N) & = & \begin{bmatrix}0 & \cdots & 0 & p_{d}(0)B^{T}S_{N+1}A\end{bmatrix}\\
A_{11}(N-1) & = & 2\bar{P}_{d}(0)P(0)B^{T}S_{N+1}AB+T_{c}(1,1)\\
A_{12}(N-1,\tau_{N-2})\zeta_{N-2} & = & M(N-1)x_{N-1}+K_{\zeta}(1,\tau_{N-2},1)
\end{eqnarray*}

For $k>1$

\begin{align*}
A_{11}(N-k)= & -K_{gg}(k)-2p(0)K_{\theta RL}(k-1,1)B\\
 &+T_{c}(k,1)+2K_{e}(k)\\
A_{12}(N-k,\tau_{N-k-1})\zeta_{N-k-1}= & K_{\zeta}(k,\tau_{N-k-1},1)-K_{gx}(k,1)x_{N-k}\\
 & -K_{gu}(k,\tau_{N-k-1},1,1)+M(N-k)x_{N-k}\\
 & 
\end{align*}

\[
\zeta_{N-k-1}=\left[\begin{array}{c}
v_{N-k-1}\\
\vdots\\
v_{\tau_{N-k-1}}\\
x_{N-k}
\end{array}\right]
\]

\subsection*{Expressions for $\boldsymbol{T}$}

\begin{eqnarray*}
T(0,b) & = & \sum_{i=0}^{N+1-b}v_{i}^{T}R_{i}v_{i}+x_{N}^{T}A^{T}S_{N+1}Ax_{N}+\sum_{i=0}^{N+1-b}p_{d}(N-i)v_{i}^{T}B^{T}S_{N+1}Bv_{i}\\
R_{i} & = & P_{d}(N-i)R\\
T(k,b) & = & T_{noX}(k,b)+x_{N-k}^{T}T_{X}(k,b)x_{N-k}\\
T(k,b) & = & T(k,b+1)+v_{N+1-b}^{T}T_{c}(k,b)v_{N+1-b}\\
T_{c}(0,b) & = & R_{N+1-b}+p_{d}(b-1)B^{T}S_{N+1}B\\
T_{noX}(0,b) & = & \sum_{i=0}^{N+1-b}v_{i}^{T}R_{i}v_{i}+\sum_{i=0}^{N+1-b}p_{d}(N-i)v_{i}^{T}B^{T}S_{N+1}Bv_{i}\\
T_{X}(0) & = & A^{T}S_{N+1}A
\end{eqnarray*}

For $k>0$

\begin{eqnarray*}
T(k,b) & = & T_{noX}(k-1,b+1)+\sum_{i=0}^{N-k+1-b}p_{d}(N-k-i)v_{i}^{T}B^{T}S_{N-k+1}Bv_{i}\\
 &  & +x_{N-k}^{T}A^{T}S_{N-k+1}Ax_{N-k}\\
T_{c}(k,b) & = & T_{c}(k-1,b+1)+p_{d}(b-1)B^{T}S_{N-k+1}B\\
T_{noX}(k,b) & = & T_{noX}(k-1,b+1)+\sum_{i=0}^{N-k+1-b}p_{d}(N-k-i)v_{i}^{T}B^{T}S_{N-k+1}Bv_{i}\\
T_{X}(k) & = & A^{T}S_{N-k+1}A
\end{eqnarray*}

\subsection*{Expressions for $\boldsymbol{S_{N}}$}
If $k\leq1$

\[
S_{N-k}=T_{X}(k)-M(N-k)^{T}A_{11}^{-1}(N-k)M(N-k)+Q
\]

If $k>1$

\begin{eqnarray*}
S_{N-k} & = & -(M(N-k)-K_{gx}(k,\tau_{N-k-1}))^{T}A_{11}^{-1}(N-k)(M(N-k)-K_{gx}(k,\tau_{N-k-1}))\\
 &  & +T_{X}(k)+Q\\
\end{eqnarray*}

\subsection*{Expressions for $\boldsymbol{M}$}

\begin{eqnarray*}
M(N-k) & = & \sum_{i=N-k}^{N}p_{d}(i-(N-k)) B^{T}S_{i+1}(A)^{i-(N-k)+1}\\
\end{eqnarray*}

\subsection*{Expressions for $\boldsymbol{K_{e}}$}

\begin{align*}
K_{e}(k)= & \sum_{i=N-k+1}^{N}\sum_{j=N-k+1}^{i}\bar{P}_{d}(i-(N-k+1))\\
 & \times P(k+j-N-1)B^{T}S_{i+1}A^{i-j+1}B
\end{align*}

\subsection*{Expressions for $\boldsymbol{K_{gg}}$}

For $k<3$

\begin{align*}
K_{gg}(k)= & p(0)K_{\eta RL}^{T}(k-1,1)A_{11}^{-1}(N-k+1)K_{\eta RL}(k-1,1)\\
 & +\bar{p}(0)K_{\eta C}^{T}(k-1,1)A_{11}^{-1}(N-k+1)K_{\eta C}(k-1,1)\\
\end{align*}

For $k\geq3$

\begin{align*}
K_{gg}(k)= & p(0)K_{\eta RL}^{T}(k-1,1)A_{11}^{-1}(N-k+1)K_{\eta RL}(k-1,1)\\
 & +\bar{p}(0)K_{\eta C}^{T}(k-1,1) A_{11}^{-1}(N-k+1)K_{\eta C}(k-1,1)\\
 & +2\bar{p}(0)K_{uuCs}(k-1,1) +\bar{p}(0)K_{rCs}(k-1,1)\\
 & +p(0)\left(2K_{uuRL}(k-1,1)+K_{rRL}(k-1,1)\right)\\
\end{align*}

\subsection*{Expressions for $\boldsymbol{K_{\alpha}}$}

\begin{align*}
K_{\alpha}(k,\tau,b)= & \sum_{j=\tau+1}^{N-k-b}p(N-k-j) K_{\eta}^{T}(k-1,j,b+1)\\
 & \times A_{11}^{-1}(N-k+1)K_{\eta}(k-1,j,b+1)\\
 & +\bar{P}(N-k-1-\tau)K_{\eta}^{T}(k-1,\tau,b+1)\\
 & \times A_{11}^{-1}(N-k+1)K_{\eta}(k-1,\tau,b+1)\\
\end{align*}
In the case when $\tau<N-k-b$ then

\begin{align*}
  K_{\alpha}(k,\tau,b)=&v_{N-k-b}^{T}K_{\alpha Cs}(k,b)v_{N-k-b}\\
 & +K_{\alpha Cd}^{T}(k,\tau,b+1,b)v_{N-k-b}\\
 & +v_{N-k-b}^{T}K_{\alpha Cd}(k,\tau,b+1,b)\\
 & +K_{\alpha}(k,\tau,b+1).
\end{align*}
The subfunctions above are given by

\begin{alignat*}{1}
 K_{\alpha Cs}(k,b)=&p(b)K_{\eta RL}^{T}(k-1,b+1)A_{11}^{-1}(N-k+1) \\
 & \times K_{\eta RL}(k-1,b+1)\\
 & +\bar{P}(b)K_{\eta C}^{T}(k-1,b+1)A_{11}^{-1}(N-k+1)\\
 & \times K_{\eta C}(k-1,b+1)\\
\\
 K_{\alpha RL}(k,b)=&\bar{P}(b-\text{1})K_{\eta RL}^{T}(k-1,b+1)\\
 & \times A_{11}^{-1}(N-k+1)K_{\eta RL}(k-1,b+1)
 \\
  K_{\alpha Cd}(k,\tau,b,h)=&\sum_{j=\tau+1}^{N-k-b}p(N-k-j)K_{\eta C}^{T}(k-1,h+1)\\
 & \times A_{11}^{-1}(N-k+1)K_{\eta}(k-1,j,b+1)\\
 & +\bar{P}(N-k-1-\tau)K_{\eta C}^{T}(k-1,h+1)\\
 & \times A_{11}^{-1}(N-k+1)K_{\eta}(k-1,\tau,b+1)
\end{alignat*}
In the case when $\tau<N-k-b$
\begin{align*}
  K_{\alpha Cd}(k,\tau,b,h)= & K_{\alpha Cd}(k,\tau,b+1,h)\\
 & +K_{\alpha CdC}(k,b,h)v_{N-k-b}
\end{align*}
Furthermore,

\begin{align*}
 K_{\alpha CdC}(k,b,h)= & \bar{P}(b)K_{\eta C}^{T}(k-1,h+1)\\
 & \times A_{11}^{-1}(N-k+1)K_{\eta C}(k-1,b+1)\\
 & +p(b)K_{\eta C}^{T}(k-1,h+1)\\
 & \times A_{11}^{-1}(N-k+1)K_{\eta RL}(k-1,b+1)\\
\\
  K_{\alpha CdRL}(k,b,h)\text{=} & \bar{P}(b-1)K_{\eta C}^{T}(k-1,h+1)\\
 & \times A_{11}^{-1}(N-k+1)K_{\eta RL}(k-1,b+1).
\end{align*}

\subsection*{Expressions for $\boldsymbol{K_{\beta}}$}

\begin{align*}
 K_{\beta}(k,\tau,b)= & \left(\sum_{j=\tau+1}^{N-k-b}p(N-k-j)K_{\eta}(k-1,j,b+1)\right.\\
 & +\bar{P}(N-k-1-\tau_{N-3})K_{\eta}(k-1,\tau,b+1)\bigg)^{T}\\
 & \times A_{11}^{-1}(N-k+1)K_{\eta C}(k-1,1)\\
 K_{\beta C}(k,b)= & \left(p(b)K_{\eta RL}(k-1,b+1)\right.\\
 & \left.+\bar{P}(b)K_{\eta C}(k-1,b+1)\right)^{T}\\
 & \times A_{11}^{-1}(N-k+1)K_{\eta C}(k-1,1)\\
  K_{\beta RL}(k,b)= & \bar{P}(b-1)K_{\eta RL}^{T}(k-1,b+1)\\
  &\times A_{11}^{-1}(N-k+ 1)K_{\eta C}(k-1,1).\\
\end{align*}
For $\tau<N-k-b$

\begin{align*}
 K_{\beta}(k,\tau,b)= & v_{N-k-b}^{T}K_{\beta C}(k,b)\\
 & +K_{\beta}(k,\tau,b+1)
\end{align*}

\subsection*{Expressions for $\boldsymbol{K_{r}}$}

\[
K_{r}(2,\tau,1)=K_{\alpha}(2,\tau,1)
\]
For $k>2$

\begin{alignat*}{1}
 K_{r}(k,\tau,b)=& \sum_{j=\tau+1}^{N-k-b}p(N-k-j)K_{r}(k-1,j,b+1)\\
 & +\bar{P}(N-k-1-\tau)K_{r}(k-1,\tau,b+1)\\
 & +K_{\alpha}(k,\tau,b)
\end{alignat*}
For $\tau<N-k-b$

\begin{alignat*}{1}
K_{r}(k,\tau,b)= & K_{r}(k,\tau,b+1)\\
 & +v_{N-k-b}^{T}K_{rCs}(k,b)v_{N-k-b}\\
 & +2v_{N-k-b}^{T}K_{rCd}(k,\tau,b+1,b)\\
\end{alignat*}
Moreover,

\begin{align*}
K_{rCs}(k,b) = & \bar{P}(b)K_{rCs}(k-1,b+1)\\ &+p(b)K_{rRL}(k-1,b+1)\\
  &    +K_{\alpha Cs}(k,b)\\
 K_{rRL}(k,b)= & \bar{P}(b-1)K_{rRL}(k-1,b+1)+K_{\alpha RL}(k,b)\\
 K_{rCd}(k,\tau,b,h)= & \sum_{j=\tau+1}^{N-k-b}p(N-k-j)\\
 &\times K_{rCd}(k-1,j,b+1,h+1)\\
 & +\bar{P}(N-k-1-\tau_{N-4})\\
 &\times K_{rCd}(k-1,\tau,b+1,h+1)\\
 & +K_{\alpha Cd}(k,\tau,b,h)
\end{align*}

For $\tau<N-k-b$

\begin{eqnarray*}
K_{rCd}(k,\tau,b,h) & = & K_{rCd}(k,\tau,b+1,h)\\
 &  & +K_{rCdC}(k,b,h)v_{N-k-b}
\end{eqnarray*}
Furthermore,

\begin{align*}
K_{rCdC}(k,b,h)= & \bar{P}(b)K_{rCdC}(k-1,b+1,h)\\
 & +K_{\alpha CdC}(k,b,h)\\
 & +p(b)K_{rCdRL}(k-1,b+1,h)\\
 K_{rCdRL}(k,b,h)= & \bar{P}(b-1)K_{rCdRL}(k-1,b+1,h)\\
 & +K_{\alpha CdRL}(k,b,h)
\end{align*}

\subsection*{Expressions for $\boldsymbol{K_{\eta}}$}

\begin{eqnarray*}
K_{\eta}(1,\tau,b) & = & K_{\zeta}(1,\tau,b)\\
K_{\eta C}(1,b) & = & K_{\zeta C}(1,b)\\
K_{\eta RL}(1,b) & = & K_{\zeta RL}(1,b)
\end{eqnarray*}
For $k>1$

\begin{eqnarray*}
K_{\eta}(k,\tau,b) & = & K_{\zeta}(k,\tau,b)\\
 &  & -K_{gu}(k,\tau,b,1)\\
K_{\eta C}(k,b) & = & K_{\zeta C}(k,b)-K_{guC}(k,b)\\
K_{\eta RL}(k,b) & = & K_{\zeta RL}(k,b)-K_{guRL}(k,b)
\end{eqnarray*}
For $\tau<N-k-b$

\begin{eqnarray*}
K_{\eta}(k,\tau,b) & = & K_{\eta}(k,\tau,b+1)\\
 &  & +K_{\eta C}(k,b)v_{N-k-b}
\end{eqnarray*}

\subsection*{$\boldsymbol{K_{\zeta}}$}

\begin{align*}
K_{\zeta}(k,\tau,b)=& \sum_{i=N-k+1}^{N}\sum_{j=N-k+1}^{i}\bar{P}_{d}(i-(N-k+1))B^{T}S_{i+1}A^{i-j+1}\\
 &   \times B(P(i-N+k)-P(k+j-N-1))\\
 &   \times\left(\sum_{t=b}^{N-k-1-\tau}\frac{P(k+j-N-1+t)-P(t-1)}{\bar{P}(t-1)}\right.\\
 &   \times\prod_{h=2}^{t}(\frac{\bar{P}(k+j-N+h-2)}{\bar{P}(h-2)})v_{N-k-t}\\
 &   \left.+\prod_{h=1}^{N-k-1-\tau}(\frac{\bar{P}(k+j-N+h-1)}{\bar{P}(h-1)})v_{\tau}\right)
\end{align*}
For $\tau<N-k-b$

\begin{align*}
  K_{\zeta}(k,\tau,b)= &K_{\zeta}(k,\tau,b+1)+K_{\zeta C}(k,b)v_{N-k-b}\\
 K_{\zeta C}(k,b)= &  \sum_{i=N-k+1}^{N}\sum_{j=N-k+1}^{i}\bar{P}_{d}(i-(N-k+1))B^{T}S_{i+1}\\
 &  \times A^{i-j+1}B(P(i-N+k)-P(k+j-N-1))\\
 &  \times\left(\frac{P(k+j-N-1+b)-P(b-1)}{\bar{P}(b-1)}\right.\\
 & \left. \times\prod_{h=2}^{b}(\frac{\bar{P}(k+j-N+h-2)}{\bar{P}(h-2)})v_{N-k-b}\right)\\
K_{\zeta RL}(k,b)=& \sum_{i=N-k+1}^{N}\sum_{j=N-k+1}^{i}\bar{P}_{d}(i-(N-k+1))B^{T}S_{i+1}\\
 &   \times A^{i-j+1}B(P(i-N+k)-P(k+j-N-1))\\
 &   \times\prod_{h=1}^{b-1}(\frac{\bar{P}(k+j-N+h-1)}{\bar{P}(h-1)})
\end{align*}

\subsection*{Expressions for $\boldsymbol{K_{\theta}}$}

For $k=1$

\begin{eqnarray*}
K_{\theta}(1,\tau,b) & = & K_{\eta}^{T}(1,\tau,b)A_{11}^{-1}(N-1)M(N-1)\\
K_{\theta C}(1,b) & = & K_{\eta C}^{T}(1,b)A_{11}^{-1}(N-1)M(N-1)\\
K_{\theta RL}(1,b) & = & K_{\eta RL}^{T}(1,b)A_{11}^{-1}(N-1)M(N-1)
\end{eqnarray*}
If $\tau<N-k-b$

\begin{eqnarray*}
K_{\theta}(1,\tau,b) & = & K_{\theta}(1,\tau,b+1)\\
 &  & +v_{N-k-b}^{T}K_{\theta C}(1,b)
\end{eqnarray*}
For $k>1$

\begin{eqnarray*}
K_{\theta}(k,\tau,b) & = & K_{\eta}^{T}(k,\tau,b)A_{11}^{-1}(N-k)\\
 &  & \times\left(-K_{gx}(k,1)+M(N-k)\right)\\
 &  & +K_{ux}(k,\tau,b)\\
K_{\theta C}(k,b) & = & K_{\eta C}^{T}(k,b)A_{11}^{-1}(N-k)\\
 &  & \times\left((-K_{gx}(k,1)+M(N-k)\right)\\
 &  & +K_{uxC}(k,b)\\
K_{\theta RL}(k,b) & = & K_{\eta RL}^{T}(k,b)A_{11}^{-1}(N-k)\\
 &  & \times\left((-K_{gx}(k,1)+M(N-2)\right)\\
 &  & +K_{uxRL}(k,b)
\end{eqnarray*}
If $\tau<N-k-b$

\begin{eqnarray*}
K_{\theta}(k,\tau,b) & = & K_{\theta}(k,\tau,b+1)\\
 &  & +v_{N-k-b}^{T}K_{\theta C}(k,b)
\end{eqnarray*}

\subsection*{Expressions for $\boldsymbol{K_{ux}}$}

\begin{align*}
 K_{ux}(k,\tau,b)= & \sum_{j=\tau+1}^{N-k-b}p(N-k-j)K_{\theta}(k-1,j,b+1)A\\
 & +\bar{P}(N-k-1-\tau_{N-3})K_{\theta}(k-1,\tau,b+1)A\\
 K_{uxC}(k,b)=  & \bar{P}(b)K_{\theta C}(k-1,b+1)A\\
 &   +p(b)K_{\theta RL}(k-1,b+1)A\\
 K_{uxC}(k,b)= &  \bar{P}(b-1)K_{\theta RL}(k-1,b+1)A\
\end{align*}

If $\tau<N-k-b$

\begin{eqnarray*}
K_{ux}(k,\tau,b) & = & K_{ux}(k,\tau,b+1)\\
 &  & +v_{N-k-b}^{T}K_{uxC}(k,b)
\end{eqnarray*}

\subsection*{Expressions for $\boldsymbol{K_{gx}}$}

\begin{align*}
K_{gx}(k,b)= & \bar{p}(0)K_{\theta C}(k-1,b)A+p(0)K_{\theta RL}(k-1,b)A\\
\end{align*}

\subsection*{Expressions for $\boldsymbol{K_{uu}}$}

\begin{align*}
K_{uu}(2,\tau,b)= & \sum_{j=\tau+1}^{N-2-b}p(N-2-j)\\
 & \times K_{\theta}(1,j,b+1)Bv_{j}\\
 & +\bar{P}(N-3-\tau)\\
 & \times K_{\theta}(1,\tau,b+1)Bv_{\tau}\\
K_{uuCs}(2,\tau,b)= & p(b)K_{\theta RL}(1,b+1)B\\
K_{uuRL}(2,b)= & \bar{P}(b-1)K_{\theta RL}(1,b+1)B\\
K_{uuCd}(2,\tau,b,h)= & \sum_{j=\tau+1}^{N-2-b}p(N-2-j)\\
 & \times K_{\theta C}(1,h)Bv_{j}\\
 & +\bar{P}(N-3-\tau)\\
 & \times K_{\theta C}(1,h)Bv_{\tau}\\
K_{uuCdC}(2,\tau,b,h)= & p(b)K_{\theta C}(1,h)B\\
K_{uuCdRL}(2,b,h)= & \bar{P}(b-1)K_{\theta C}(1,h)B
\end{align*}
If $\tau<N-k-b$

\begin{align*}
 K_{uuCd}(2,\tau,b,h)= & K_{uuCd}(2,\tau,b+1,h)\\
 & +K_{uuCdC}(2,\tau,b,h)v_{N-k-b}\\
 K_{uu}(2,\tau,b)= & K_{uu}(2,\tau,b+1)\\
 & +v_{N-k-b}^{T}K_{uuCd}(2,\tau,b+1,b+1)\\
 & +v_{N-k-b}^{T}K_{uuCs}(2,b)v_{N-k-b}
\end{align*}
For $k>2$

\begin{align*}
 K_{uu}(k,\tau,b)= & \sum_{j=\tau+1}^{N-k-b}p(N-k-j)K_{\theta}(k-1,j,b+1)Bv_{j}\\
 & +\bar{P}(N-k-1-\tau)\\
 &\times K_{\theta}(k-1,\tau,b+1)Bv_{\tau}\\
 & +\sum_{j=\tau+1}^{N-k-b}p(N-k-j)K_{uu}(k-1,j,b+1)\\
 & +\bar{P}(N-k-1-\tau)K_{uu}(k-1,\tau,b+1)\\
 K_{uuRL}(k,b)=& \bar{P}(b-1)K_{\theta RL}(k-1,b+1)B\\
 & +\bar{P}(b-1)K_{uuRL}(k-1,b+1)
\end{align*}
If $\tau<N-k-b$

\begin{align*}
K_{uu}(k,\tau,b)= & K_{uu}(k,\tau,b+1)\\
 & +v_{N-k-b}^{T}K_{uuCs}(k,b)v_{N-k-b}\\
 & +v_{N-k-b}^{T}K_{uuCd}(k,\tau,b+1,b+1)\\
 K_{uuCdC}(k,b,h)= & p(b)K_{\theta C}(k-1,h)B\\
 & +\bar{P}(b)K_{uuCdC}(k-1,b+1,b+1)\\
 & +p(b)K_{uuCdRL}(k-1,b+1,b+1)
\end{align*}
Moreover,

\begin{align*}
  K_{uuCs}(k,b)= & p(b)K_{\theta RL}(k-1,b+1)B\\
 & +p(b)K_{uuRL}(k-1,b+1)\\
 & +\bar{P}(b)K_{uuCs}(k-1,b+1)\\
 K_{uuCd}(k,\tau,b,h)= & \sum_{j=\tau+1}^{N-k-b}p(N-k-j)K_{\theta C}(k-1,h)Bv_{j}\\
 & +\bar{P}(N-k-1-\tau)K_{\theta C}(k-1,h)Bv_{\tau}\\
 & +\sum_{j=\tau+1}^{N-k-b}p(N-k-j)\\
 & \times K_{uuCd}(k-1,j,b+1,b+1)\\
 & +\bar{P}(N-k-1-\tau)\\
 & \times K_{uuCd}(k-1,\tau,b+1,b+1)\\
  K_{uuCd}(k,\tau,b,h)= & K_{uuCd}(k,\tau,b+1,h)\\
 & +K_{uuCdC}(k,\tau,b,h)v_{N-k-b}\\
 K_{uuCdRL}(k,b,b)= & \bar{P}(b-1)K_{\theta C}(k-1,h)B\\
 & +\bar{P}(b-1)K_{uuCdRL}(k-1,b+1,b+1)
\end{align*}

\subsection*{Expressions for $\boldsymbol{K_{gu}}$}

\begin{eqnarray*}
K_{gu}(2,\tau,b,h)  = & \sum_{j=\tau_{N-3}+1}^{N-2-b}p(N-2-j)K_{\theta C}(1,h)Bv_{j}\\
 & +\bar{P}(N-3-\tau)K_{\theta C}(1,h)Bv_{\tau_{N-3}}\\
 & +K_{\beta}(2,\tau,b)\\
K_{guC}(2,b,h) = & p(b)K_{\theta C}(1,h)B+K_{\beta C}(2,b)\\
K_{guRL}(2,b,h)  = & \bar{P}(b-1)K_{\theta C}(1,h)B +K_{\beta RL}(2,b)
\end{eqnarray*}
if $\tau<N-k-b$

\begin{eqnarray*}
K_{gu}(2,\tau,b,h) = & K_{guC}(2,b,h)v_{N-k-b} +K_{gu}(2,\tau,b+1,h)\\
\end{eqnarray*}
for $k>2$

\begin{align*}
 K_{gu}(k,\tau,b,h)=& \sum_{j=\tau+1}^{N-k-b}p(N-k-j)K_{\theta C}(k-1,h)Bv_{j}\\
 & +\bar{P}(N-k-1-\tau)K_{\theta C}(k-1,h)Bv_{\tau}\\
 & +K_{\beta}(k,\tau,b)\\
 & +\bar{P}(N-k-1-\tau)\\
 & \times K_{uuCd}(k-1,\tau,b+1,h+1)\\
 & +\sum_{j=\tau+1}^{N-k-b}p(N-k-j)\\
 &\times K_{uuCd}(k-1,j,b+1,h+1)\\
 & +\sum_{j=\tau+1}^{N-k-b}p(N-k-j)K_{rCd}(k-1,j,b+1)\\
 & +\bar{P}(N-k-1-\tau)\\
 & \times K_{rCd}(k-1,\tau,b+1,h)\\
 K_{guC}(k,b,h)= & p(b)K_{\theta C}(k-1,h)B+K_{\beta C}(k-1,b)\\
 & +\bar{P}(b)K_{uuCdC}(k-1,b+1,h+1)\\
 & +\bar{P}(b)K_{rCdC}(k-1,b+1,h+1)\\
 & +p(b)K_{uuCdRL}(k-1,b+1)\\
 & +p(b)K_{rCdRL}(k-1,b+1)\\
  K_{guRL}(k,b)= & \bar{P}(b-1)K_{\theta C}(k-1,h)B+K_{\beta RL}(k,b)\\
 & +\bar{P}(b-1)K_{uuCdRL}(k-1,b+1)\\
 & +\bar{P}(b-1)K_{rCdRL}(k-1,b+1)
\end{align*}
if $\tau<N-k-b$

\begin{align*}
K_{gu}(k,\tau,b,h)= & K_{gu}(k,\tau,b+1,h)\\
 & +K_{guC}(k,b,h)v_{N-k-b}
\end{align*}
\bibliographystyle{plain}   
\bibliography{example2} 

\begin{thebibliography}{10}

\bibitem{astrompredict}
Karl~Johan {\AA}str{\"o}m.
\newblock {\em Introduction to stochastic control theory}, volume~70 of {\em
  Mathematics in science and engineering}.
\newblock Academic Press, 1970.

\bibitem{bengtsson2016lqg}
Fredrik Bengtsson, Babak Hassibi, and Torsten Wik.
\newblock {LQG} control for systems with random unbounded communication delay.
\newblock In {\em Decision and Control (CDC), 2016 IEEE 55th Conference on},
  pages 1048--1055. IEEE, 2016.

\bibitem{1184678}
Christoforos~N Hadjicostis and Rouzbeh Touri.
\newblock Feedback control utilizing packet dropping network links.
\newblock In {\em Proceedings of the 41st IEEE Conference on Decision and
  Control}, volume~2, pages 1205--1210. IEEE, 2002.

\bibitem{hespanha2007survey}
Joo~P Hespanha, Payam Naghshtabrizi, and Yonggang Xu.
\newblock A survey of recent results in networked control systems.
\newblock {\em Proceedings of the IEEE}, 95(1):138--162, 2007.

\bibitem{Imer20061429}
Orhan~C. Imer, Serdar Y\"uksel, and Tamer Basar.
\newblock Optimal control of {LTI} systems over unreliable communication links.
\newblock {\em Automatica}, 42(9):1429 -- 1439, 2006.

\bibitem{liang2016optimal}
Xiao Liang, Juanjuan Xu, and Huanshui Zhang.
\newblock Optimal control and stabilization for networked control systems with
  packet dropout and input delay.
\newblock {\em IEEE Transactions on Circuits and Systems II: Express Briefs},
  64(9):1087--1091, 2016.

\bibitem{lincoln2000optimal}
Bo~Lincoln and Bo~Bernhardsson.
\newblock Optimal control over networks with long random delays.
\newblock In {\em Proceedings of the International Symposium on Mathematical
  Theory of Networks and Systems}, volume~7. Perpignan, 2000.

\bibitem{ma2017optimal}
Xiao Ma, Qingyuan Qi, and Huanshui Zhang.
\newblock Optimal {LQ} control for {NCS}s with packet dropout and delay.
\newblock In {\em 2017 Chinese Automation Congress (CAC)}, pages 5726--5729.
  IEEE, 2017.

\bibitem{5717722}
M.~Moayedi, Y.K. Foo, and Y.C. Soh.
\newblock Networked {LQG} control over unreliable channels.
\newblock In {\em Decision and Control (CDC), 2010 49th IEEE Conference on},
  pages 5851--5856, Dec 2010.

\bibitem{4908929}
Luca Schenato.
\newblock To zero or to hold control inputs with lossy links?
\newblock {\em IEEE Transactions on Automatic Control}, 54(5):1093--1099, 2009.

\bibitem{4118476}
Luca Schenato, Bruno Sinopoli, Massimo Franceschetti, Kameshwar Poolla, and
  S~Shankar Sastry.
\newblock Foundations of control and estimation over lossy networks.
\newblock {\em Proceedings of the IEEE}, 95(1):163--187, 2007.

\bibitem{shousong2003stochastic}
Hu~Shousong and Zhu Qixin.
\newblock Stochastic optimal control and analysis of stability of networked
  control systems with long delay.
\newblock {\em Automatica}, 39(11):1877--1884, 2003.

\bibitem{wang2018optimal}
Zhuwei Wang, Lihan Liu, Chao Fang, Xiaodong Wang, Pengbo Si, and Hong Wu.
\newblock Optimal linear quadratic control for wireless sensor and actuator
  networks with random delays and packet dropouts.
\newblock {\em International Journal of Distributed Sensor Networks},
  14(6):1550147718779560, 2018.

\bibitem{xu2012stochastic}
Hao Xu, Sarangapani Jagannathan, and Frank~L Lewis.
\newblock Stochastic optimal control of unknown linear networked control system
  in the presence of random delays and packet losses.
\newblock {\em Automatica}, 48(6):1017--1030, 2012.

\bibitem{6760935}
Jen{-}te Yu and Li-Chen Fu.
\newblock A new compensation framework for {LQ} control over lossy networks.
\newblock In {\em 2013 IEEE 52nd Annual Conference on Decision and Control
  (CDC),}, pages 6610--6614, Dec 2013.

\end{thebibliography}

\end{document}